\documentclass[11pt]{amsart}
\usepackage{amsthm,amsmath,amssymb}
\usepackage{mathtools}
\usepackage{graphicx}
\usepackage{caption}

\graphicspath{{Figure/}}
\numberwithin{equation}{section}
\numberwithin{figure}{section}
\theoremstyle{plain}
\newtheorem{thm}{Theorem}[section]
\newtheorem{lem}[thm]{Lemma}

\newtheorem{prop}[thm]{Proposition}

\theoremstyle{remark}

\newcommand{\RomanNumeralCaps}[1]
    {\MakeUppercase{\romannumeral #1}}
\begin{document}
\setlength{\abovedisplayskip}{10pt}
\setlength{\belowdisplayskip}{10pt}

\title{Three product formulas for ratios of tiling counts of hexagons with collinear holes}

\author{Seok Hyun Byun}
\address{Department of Mahtematics, Indiana University, Bloomington}
\email{byunse@indiana.edu}

\maketitle
\begin{abstract}
Rosengren found an explicit formula for a certain weighted enumeration of lozenge tilings of a hexagon with an arbitrary triangular hole. He pointed out that a certain ratio corresponding to two such regions has a nice product formula. In this paper, we generalize this to hexagons with arbitrary collinear holes. It turns out that, by using same approach, we can also generalize Ciucu's work on the number and the number of centrally symmetric tilings of a hexagon with a fern removed from its center. This proves a recent conjecture of Ciucu.
\end{abstract}

\section{Introduction}
Enumeration of lozenge tilings of a region on a triangular lattice has been studied for many decades. In particular, people are interested in regions whose number of lozenge tilings is expressed as a simple product formula. One such region is a hexagonal region with a triangular hole in the center. Many works have been done on this topic by Ciucu [2], Ciucu et al. [6], and Okada and Krattenthaler [10]. Later, Rosengren [11] found a formula for a weighted enumeration of lozenge tilings of a hexagon with an arbitrary triangular hole. He pointed out that the ratio between numbers of lozenge tilings of two such regions whose holes have symmetric position with respect to the center has a nice product formula. In this paper, we give a conceptual explanation of the symmetry, which enables us to generalize the result to hexagons with arbitrary collinear triangular holes. In his paper, Ciucu [3] defined a new structure, called a fern, which is an arbitrary string of triangles of alternating orientations that touch at corners and are lined up along a common axis. He considered a hexagon with a fern removed from its center and proved that the ratio of the number of lozenge tilings of two such regions is given by a simple prodcut formula. Later, Ciucu [5] also proved that the same kind of ratio for centrally symmetric lozenge tilings also has a simple product formula. In particular, he pointed out that for hexagons with a fern removed from the center, the ratio of centrally symmetric lozenge tilings is the square root of the ratio of the total number of tilings. Ciucu also conjectured in [5] (See also [4]) that this square root phenomenon holds more generally, when any finite number of collinear ferns are removed in a centrally symmetric way. In this current paper, we prove Ciucu's conjecture, and we extend it further.

\section{Statement of Main Results}
Any hexagon on a triangular lattice has a property that difference between two parallel sides is equal for all 3 pairs. Thus, we can assume that the side lengths of the hexagon are \textit{a, b+k, c, a+k, b, c+k} in clockwise order, where \textit{a} is a length of a top side. Also, without loss of generality, we can assume that \textit{k} is non-negative and a southeastern side of a hexagon (=\textit{c}) is longer than or equal to a side length of a southwestern side (=\textit{b}). Note that this hexagonal region has \textit{k} more up-pointing unit-triangles than down-pointing unit-triangles. Since every lozenge consists of one up-pointing unit-triangle and one down-pointing unit-triangle, to be completely tiled by lozenges, we have to remove \textit{k} more up-pointing unit-triangles than down-pointing unit-triangles from the hexagon. There are many ways to do that, but let's consider a following case.
Let's call a set of triangles on a triangular lattice is \textit{collinear} or \textit{lined up} if horizontal side of all triangles are on a same line.
Now, let's consider any horizontal line passing through the hexagon. Suppose the line is \textit{l}-th horizontal line from bottom side of the hexagon. Note that the length of the horizontal line depends on the size of \textit{l}: Let's denote the length of the line by \textit{L(l)}. Then we have $L(l)=a+k-l+min(b,l)+min(c,l)$.

For any subsets $X=\{x_1, ..., x_{m+k}\}$ and $Y=\{y_1, ..., y_m\}$ of $[L(l)]:=\{1,2,...,L(l)\}$, let $H_{a,b,c}^{k,l}(X:Y)$ be the region obtained from the hexagon of side length $a$, $b+k$, $c$, $a+k$, $b$, $c+k$ in clockwise order from top by removing up-pointing unit-triangles whose labels of horizontal sides form a set $X=\{x_1, x_2, ..., x_{m+k}\}$, and down-pointing unit-triangles whose labels of horizontal sides form a set $Y=\{y_1, y_2, ..., y_m\}$ on the \textit{l}-th horizontal line from the bottom, where labeling on the horizontal line is $1,2,...,\textit{L(l)}$ from \textbf{left to right}. Let's call the horizontal line as a \textit{baseline} of removed triangles. Similarly, let $\overline{H}_{a,b,c}^{k,l}(X:Y)$ be a same kind of region, except that labeling on the horizontal line is $1$, $2$, ..., $L(l)$ from \textbf{right to left}. Also, for any region \textit{R} on a triangular lattice, let \textit{M(R)} be a number of lozenge tilings of the region. First theorem expresses a ratio of numbers of lozenge tilings of two such region as a simple product formula.

\begin{thm}
Let a, b, c, k, l, m be any non-negative integers such that $b \leq c$, $0 \leq l \leq b+c$ and $m \leq min(b, l, b+c-l)$. Also, let $X=\{x_1, x_2, ..., x_{m+k}\}$ and $Y=\{y_1, y_2, ..., y_m\}$ be subsets of $\lbrack L(l)\rbrack=\{1, 2, ..., L(l)\}$. Then
\begin{equation}
\begin{aligned}
    &\frac{M({H_{a,b,c}^{k,l}(X:Y)})}{M({\overline{H}_{a,b,c}^{k,b+c-l}(X:Y)})}\\
    &=\frac{H(k+l)H(b+c-l)}{H(l)H(b+c+k-l)}\\
    &\cdot\frac{\prod_{i=1}^{m+k}(x_i-b+max(b, l))_{(b-l)}\cdot(a+k+min(b, l)+1-x_i)_{(c-l)}}{\prod_{j=1}^{m}(y_i-b+max(b, l))_{(b-l)}\cdot(a+k+min(b, l)+1-y_j)_{(c-l)}}
\end{aligned}
\end{equation}
\end{thm}

where the hyperfactorial \textit{H}(\textit{n}) is defined by
\begin{equation}
    H(n):=0!1!\cdot\cdot\cdot(n-1)!
\end{equation}

\begin{figure}
    \centering
    \includegraphics[width=1\textwidth]{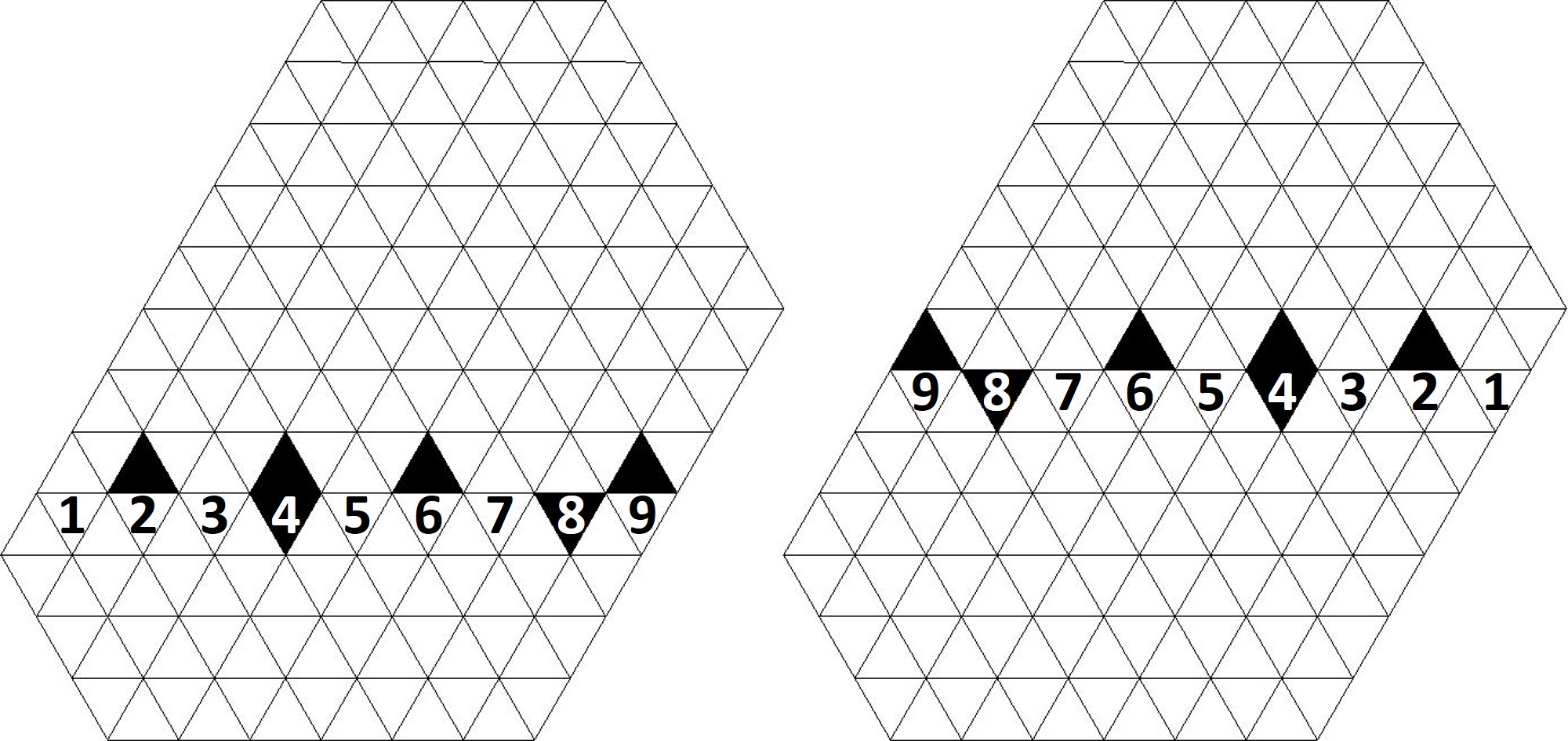}
    \caption{Two regions $H_{4,3,7}^{2,4}(\{2,4,6,9\}:\{4,8\})$(left) and $H_{4,3,7}^{2,6}(\{2,4,6,9\}:\{4,8\})$(right)}
\end{figure}

To state next results, we need to recall a result of Cohn, Larsen and Propp [8], which is a lozenge tilings interpretation of a classical result of Gelfand and Tsetlin [9]. Recall that $\Delta(S):=\prod_{s_1<s_2, s_1,s_2\in S}{(s_2-s_1)}$ and $\Delta(S,T):=\prod_{s \in S ,t \in T}{|t-s|}$ for any finite sets S and T.

\begin{prop}
For any non-negative integers $m, n$ and any subset $S=\{s_1, s_2,...,s_n\} \subset [m+n]:=\{1, 2,..., m+n\}$, let $T_{m,n}(S)$ be the region on a triangular lattice obtained from the trapezoid of side lengths $m$, $n$, $m+n$, $n$ clockwise from the top by removing the up-pointing unit-triangles whose bottoms sides are labeled by elements of a set $S=\{s_1, s_2,...,s_n\}$, where bottom side of the trapezoid is labeled by $1, 2, ..., m+n$ from left to right. Then
\begin{equation}
    M(T_{m,n}(S))=\frac{\Delta(S)}{\Delta([n])}=\frac{\Delta(S)}{H(n)}
\end{equation}
\end{prop}

For any finite subset of integers $S=\{s_1, s_2,..., s_n\}$, where elements are written in increasing order, let $T(S)$ be a region obtained by translating a region $T_{s_n-(s_1-1)-n,n}(s_1-(s_1-1), s_2-(s_1-1),...,s_n-(s_1-1))$ by $(s_1-1)$ units to the right and $s(S):=M(T(S))=\frac{\Delta(S)}{H(n)}$. A region on a triangular lattice is called \textit{balanced} if it contains same number of up-pointing unit-triangles and down-pointing unit-triangles. Geometrically, $T(S)$ is the balanced region that can be obtained from a trapezoid of bottom length $(s_n-s_1+1)$ by deleting up-pointing unit-triangles whose labels are $s_1$, $s_2$,..., $s_n$ on bottom, where bottom line is labeled by $s_1$, $(s_1+1)$,..., $(s_n-1)$, $s_n$ (See Figure 2.2).

\begin{figure}
    \centering
    \includegraphics[width=1\textwidth]{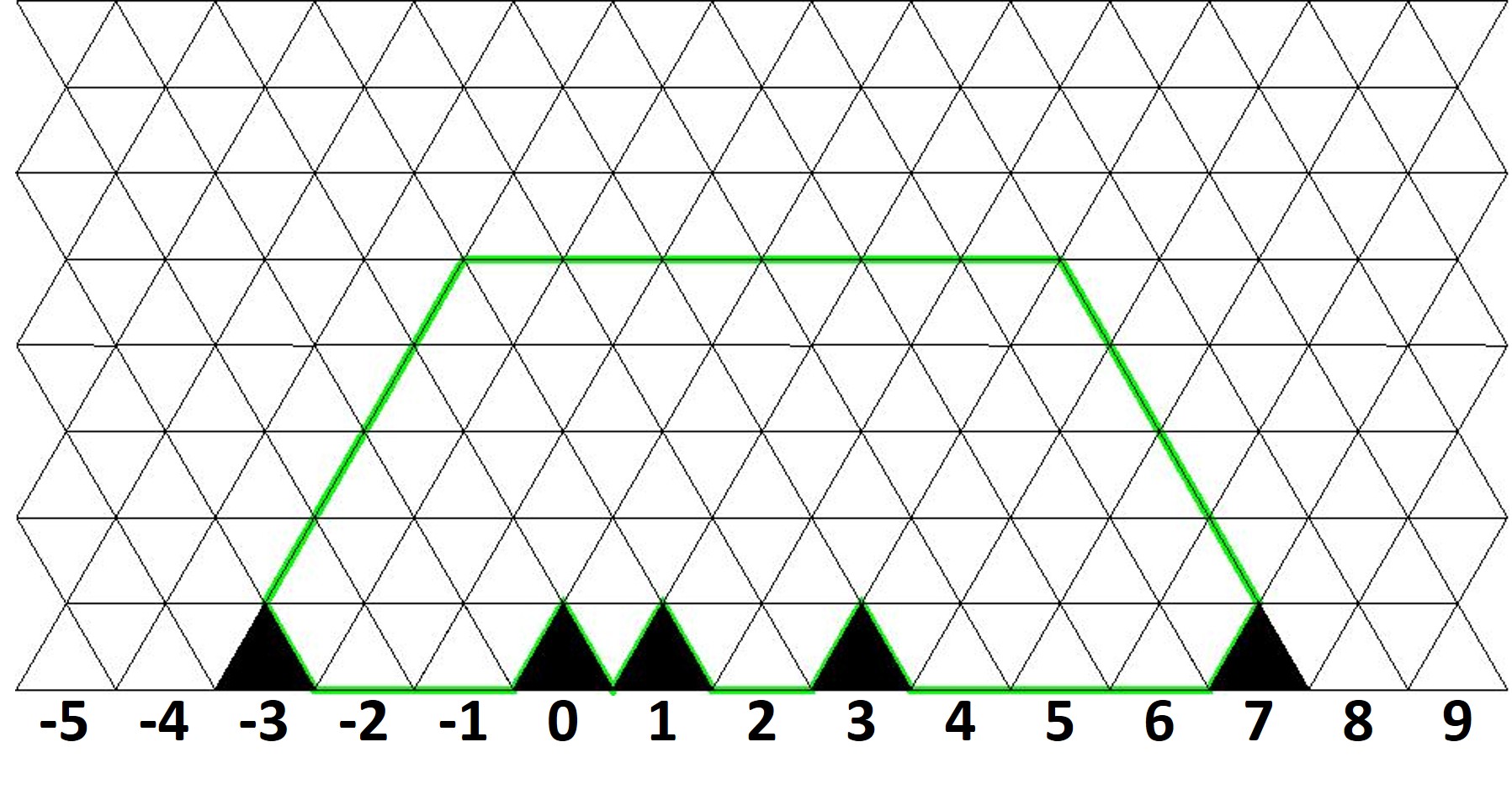}
    \caption{A region T(\{-3, 0, 1, 3, 7\})}
\end{figure}

In his paper, Ciucu [3] defined a new structure, called a \textit{fern}, which is an arbitrary string of triangles of alternating orientations that touch at corners and are lined up. For non-negative integers $a_1$,...,$a_k$, a \textit{fern} $F(a_1,...,a_k)$ is a string of \textit{k} lattice triangles lined up along a horizontal lattice line, touching at their vertices, alternately oriented up and down and having sizes $a_1$,...,$a_k$ from left to right (with the leftmost oriented up). We call the horizontal lattice line as a \textit{baseline} of the fern.

\begin{figure}
    \centering
    \includegraphics[width=1\textwidth]{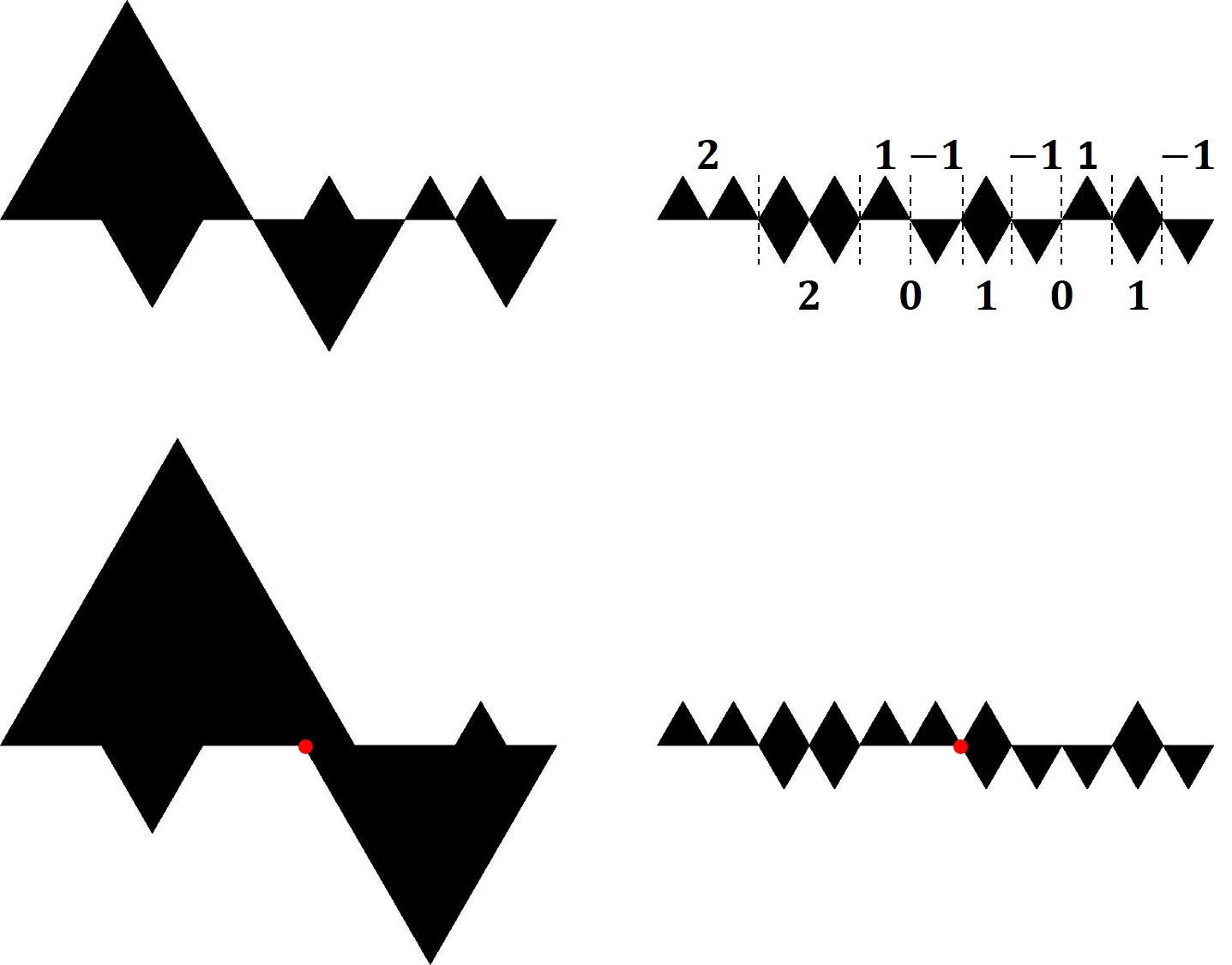}
    \caption{Budded fern $F(2,1,-1,-1,1,-1:2,0,1,0,1)$ (top left) and its baseline representation (top right), corresponding budded bowtie (bottom left) and its baseline representation (bottom right). A red point is a turning point}
    \label{fig:my_label}
\end{figure}

Now, let's give additional structure to fern by adding \textit{buds} (triangles) on the baseline, and we will call this new structure as a \textit{budded fern}. To label this new structure, we remove all unit-triangles from the budded fern, excepts unit-triangles whose horizontal side is on the baseline. We call it as a \textit{baseline representation} of the budded fern. Then we count numbers of consecutive up-pointing unit-triangles, down-pointing unit-triangles and vertical unit-lozenges on the baseline. When we count these numbers, up-pointing (or down-pointing) unit-triangle which is contained in a vertical unit-lozenge is not considered as an up-pointing (or down-pointing) unit-triangle. If an up-pointing unit-triangle and a down-pointing unit-triangle are adjacent, we think as if there are $0$ vertical lozenges between them. Now, we line up these numbers from left to right, and put - to numbers that represent numbers of down-pointing unit-triangles. Then, by allowing $a^k_1=0$ and $a^k_{r_k}=0$, we get a sequence of integers $a^k_1$, $w^k_1$, $a^k_2$, $w^k_2$,..., $a^k_{r_k-1}$, $w^k_{r_k-1}$, $a^k_{r_k}$, where $a^k_i$ represent a (signed) number of consecutive up-pointing (or down-pointing) unit-triangles, and $w^k_i$ represent a number of consecutive vertical unit-lozenges. Let $A_k$ represents a sequence $(a^k_1, a^k_2,,..., a^k_{r_k})$ and $W_k$ represents a sequence $(w^k_1, w^k_2,,..., w^k_{r_k-1})$. Then we denote the original budded fern as $F(A_k:W_k)$, and its baseline representation as $F_{br}(A_k:W_k)$. Let $L^k_i$ be a leftmost vertex of $a^k_i$ consecutive triangles, and $R^k_i$ be a rightmost vertex of the consecutive triangles. Also, let $I^k:=\{i\in[r_k]|a^k_i> 0\}$, $J^k:=\{i\in[r_k]|a^k_i< 0\}$, $p_k:=\sum_{i\in I^k} a^k_i$ and $n_k:=-\sum_{i\in J^k} a^k_i$.
From this budded fern, we will construct a corresponding \textit{budded bowtie} as follows. From a baseline representation of the budded fern, we move up-pointing unit-triangles to left, and down-pointing unit-triangles to right along the baseline, fixing vertical lozenges. Then we call a right vertex of a right-most up-poiting unit-triangle which is not a part of a vertical lozenge as a \textit{turning point} of a new structure and denote it by $T^k$. Then we put vertical lozenges between consecutive up-pointing (or down-pointing) unit-triangles as much as possible. Then we get a bowtie (possibly a slipped bowtie) with some triangles attached. We call this as a \textit{budded bowtie} and denote it and its baseline representation by $B(A_k:W_k)$ and $B_{br}(A_k:W_k)$, respectively. Also, let $u_k$ be a smallest positive integer such that $p_k \leq |a^k_1|+|a^k_2|+...+|a^k_{u_k}|$, and $v_k \in [a^k_{u_k}]$ be a positive integer such that $|a^k_1|+|a^k_2|+...+|a^k_{u_k-1}|+v_k=p_k$.

When we say a budded fern $F(A_k:W_k)$, we equip with corresponding sequences $A_k$, $W_k$, sets $I^k$, $J^k$, indices $r_k$, $p_k$, $n_k$, $u_k$, $v_k$ and vertices $L_1^k, L_2^k,..., L_{r_k}^k$, $R_1^k, R_2^k,..., R_{r_k}^k$, $T^k$.

Now, let $H_{a,b,c}^{k,l}(F(A_1:W_1),..., F(A_t:W_t) : m_1, m_2,..., m_t, m_{t+1})$ be a region obtained from the hexagon of side length $a$, $b+k$, $c$, $a+k$, $b$, $c+k$ in clockwise order from top by removing budded ferns $F(A_1:W_1),..., F(A_n:W_n)$ on a $l$-th horizontal line from the bottom so that a distance between a leftmost vertex on the horizontal line and a leftmost vertex of $F(A_1:W_1)$ is $m_1$, a distance between a rightmost vertex on the horizontal line and a rightmost vertex of $F(A_t:W_t)$ is $m_{t+1}$, and a distance between two adjacent budded ferns $F(A_i:W_i)$ and $F(A_{i+1}:W_{i+1})$ is $m_{i+1}$ for all $i\in[t-1]$. We can similarly think of a region $H_{a,b,c}^{k,l}(F_{br}(A_1:W_1),..., F_{br}(A_t:W_t) : m_1, m_2,..., m_t, m_{t+1})$. From the region, we label the $l$-th horizontal line from the bottom by 1,2,...,$L(l)$ from left to right. Let $X^1$ be a set of labels whose corresponding segment is a side of an up-pointing unit-triangular hole, but not a side of an down-pointing unit-triangular hole. Similarly, let $X^2$ be a set of labels of segments whose corresponding segment is a side of an down-pointing unit-triangular hole, but not a side of an up-pointing unit-triangular hole and $W$ be a set of label of segments whose corresponding segment is a side of both up-pointing and down-pointing unit-triangular holes.

Similarly, let $H_{a,b,c}^{k,l}(B(A_1:W_1),..., B(A_t:W_t) : m_1, m_2,..., m_t, m_{t+1})$ be a region obtained from the hexagon of side length $a$, $b+k$, $c$, $a+k$, $b$, $c+k$ in clockwise order from top by removing budded bowties $B(A_1:W_1),..., B(A_t:W_t)$ from a $l$-th horizontal line, where positions of removed budded bowties on the horizontal line is exactly same as positions of corresponding budded ferns. Again, we can think of a region $H_{a,b,c}^{k,l}(B_{br}(A_1:W_1),..., B_{br}(A_t:W_t) : m_1, m_2,..., m_t, m_{t+1})$ and we can define sets $Y^1$ and $Y^2$ from this region as we defined sets $X^1$ and $X^2$ from $H_{a,b,c}^{k,l}(F_{br}(A_1:W_1),..., F_{br}(A_t:W_t) : m_1, m_2,..., m_t, m_{t+1})$. Note that we have $X^1\cup X^2=Y^1\cup Y^2$.

For any point and a line on the triangular lattice, let distance between a point and a line be a shortest length of a path from the point to the extension of the line along lattice. Especially, for any lattice point $E$ in a hexagon, let $d_{NW}(E)$ be a distance between the point $E$ and a northwestern side of the hexagon. Similarly, we can define $d_{SW}(E)$, $d_{NE}(E)$ and $d_{SE}(E)$ to be distances between a point $E$ and a southwestern side, northeastern side and southeastern side of the hexagon, respectively. Next theorem expresses a ratio of numbers of lozenge tilings of the two regions as a simple product formula.

\begin{figure}
    \centering
    \includegraphics[width=1\textwidth]{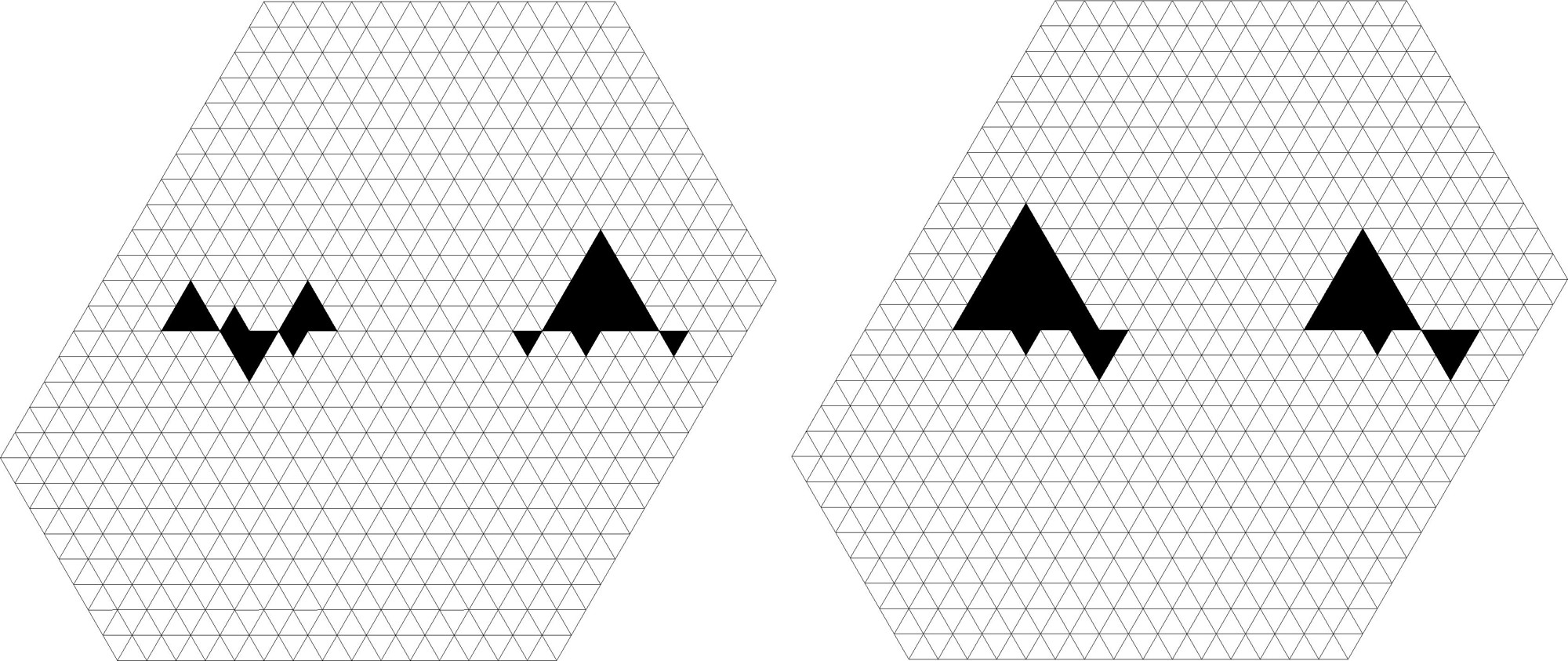}
    \caption{An example of regions (from left to right): $H_{12,8,15}^{3,13}(F(2,-1,1:1,1), F(-1,1,2,-1:0,1,0):3,6,2)$ and $H_{12,8,15}^{3,13}(B(2,-1,1 : 1,1), B(-1,1,2,-1 : 0,1,0) : 3, 6, 2)$}
\end{figure}

\begin{thm}
Let $a, b, c, k, l, m_1,...,m_{t+1}$ be any non-negative integers and $F(A_1:W_1)$,..., $F(A_t:W_t)$ be any budded ferns. Let $p:=\sum_{i=1}^{t}{p_i}$, $n:=\sum_{i=1}^{t}{n_i}$, $w:=\sum_{i=1}^{t}{\sum_{j=1}^{r_i-1}w_j^i}$ and $m:=\sum_{i=1}^{t+1}m_i$. Suppose indices satisfy following conditions: 1) $p=n+k$, 2) $p+n+w+m=L(l)$, 3) $n+w\leq min(b,l,b+c-l)$. Then we have

\begin{equation}
\begin{aligned}
    &\frac{M(H_{a,b,c}^{k,l}(F(A_1:W_1), F(A_2:W_2),..., F(A_t:W_t) : m_1, m_2,..., m_t, m_{t+1}))}{M(H_{a,b,c}^{k,l}(B(A_1:W_1), B(A_2:W_2),..., B(A_t:W_t) : m_1, m_2,..., m_t, m_{t+1}))}\\
    &=\frac{s(X^1)s(X^2)}{s(Y^1)s(Y^2)}\\
    &
\begin{aligned}
    \cdot\prod_{i=1}^{t}\Bigg[&\frac{H(d_{SW}(T^i))H(d_{NW}(L^i_{u_i}))H(d_{SE}(L^i_{u_i}))H(d_{NE}(T^i))}{H(d_{SW}(L^i_{u_i}))H(d_{NW}(T^i))H(d_{SE}(T^i))H(d_{NE}(L^i_{u_i}))}\\
    &\cdot\prod_{j < u_i, j \in J_i}\frac{H(d_{SW}(R^i_j))H(d_{NW}(L^i_j))H(d_{SE}(L^i_j))H(d_{NE}(R^i_j))}{H(d_{SW}(L^i_j))H(d_{NW}(R^i_j))H(d_{SE}(R^i_j))H(d_{NE}(L^i_j))} \\
    &\cdot\prod_{j \geq u_i, j \in I_i}\frac{H(d_{SW}(L^i_j))H(d_{NW}(R^i_j))H(d_{SE}(R^i_j))H(d_{NE}(L^i_j))}{H(d_{SW}(R^i_j))H(d_{NW}(L^i_j))H(d_{SE}(L^i_j))H(d_{NE}(R^i_j))} \Bigg]    
\end{aligned}
\end{aligned}
\end{equation}
\end{thm}

Now, let's consider a case when  a region $H_{a,b,c,k,l}(F(A_1:W_1),..., F(A_t:W_t) : m_1, m_2,..., m_t, m_{t+1})$ is centrally symmetric, that is invariant under $180^{\circ}$ rotation with respect to a center of a hexagon. To satisfy this condition, the region should satisfy following conditions:\\
(1) $b$ and $c$ have same parity\\
(2) $k=0$ and $l=\frac{b+c}{2}$\\
(3) $m_s=m_{t+2-s}$ for all $s\in [t+1]$ and $r_i=r_{t+1-i}$ for all $i\in [t]$\\
(4) $a^i_{j}=-a^{t+1-i}_{r_i+1-j}$ and $w^i_u=w^{t+1-i}_{r_i-u}$ for all $i\in [t]$, $j\in [r_i]$, $u\in [r_i-1]$

When these conditions hold, a region $H_{a,b,c}^{0,\frac{b+c}{2}}(F(A_1:W_1),..., F(A_t:W_t) : m_1, m_2,..., m_t, m_{t+1})$ and a corresponding region $H_{a,b,c}^{0,\frac{b+c}{2}}(B(A_1:W_1),..., B(A_t:W_t) : m_1, m_2,..., m_t, m_{t+1})$ are centrally symmetric, so we can compare their number of centrally symmetric lozenge tilings. Let $M_\odot(G)$ be a number of centrally symmetric lozenge tiling of a region G on a triangular lattice. The last theorem expresses a ratio of numbers of centrally symmetric lozenge tilings of the two regions as a simple product formula.

\begin{figure}
    \centering
    \includegraphics[width=1\textwidth]{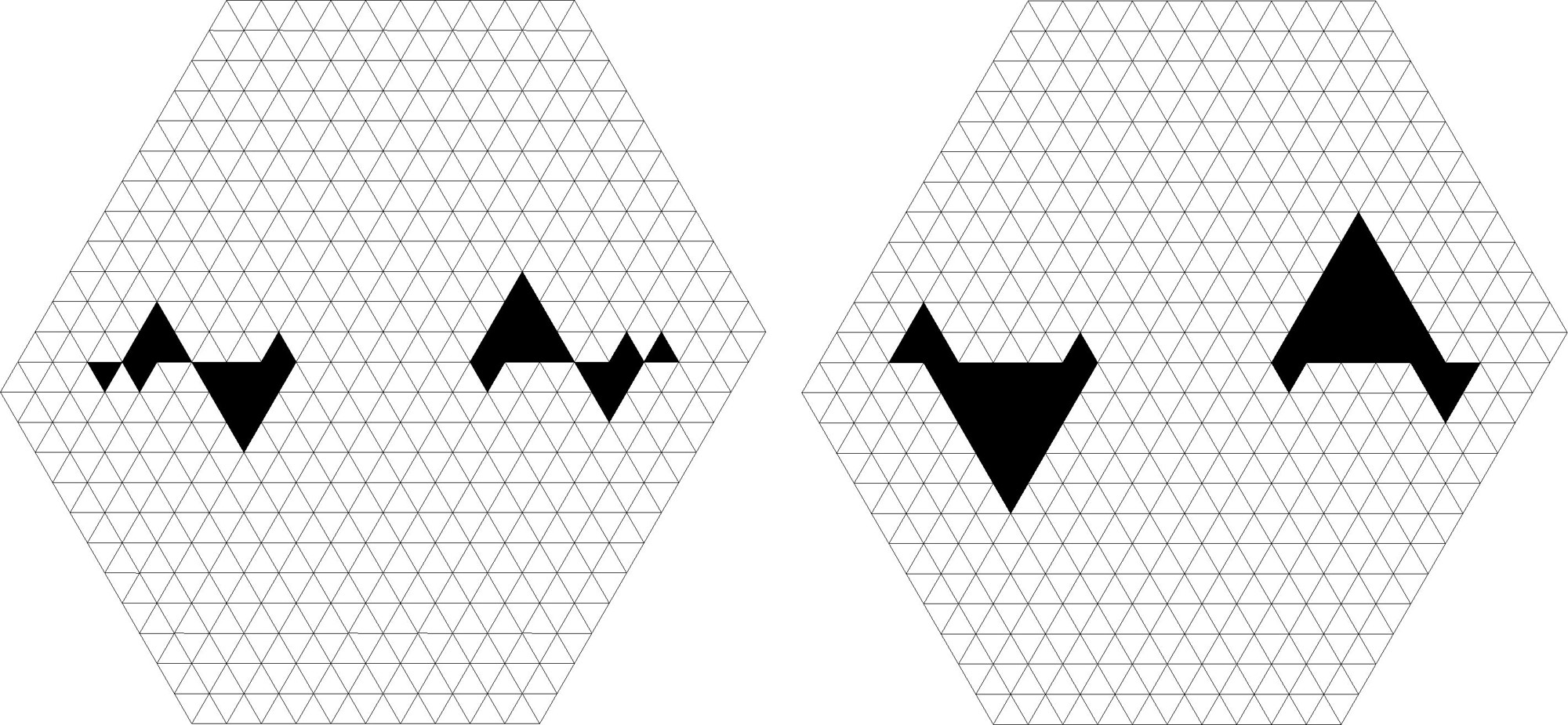}
    \caption{Centrally symmetric regions (from left to right) $H_{10,11,13}^{0,12}(F(-1,1,-2,0:1,0,1), F(0,2,-1,1:1,0,1):2,5,2)$ and $H_{10,11,13}^{0,12}(B(-1,1,-2,0:1,0,1), B(0,2,-1,1:1,0,1):2,5,2)$}
\end{figure}

\begin{thm}
Let $a,b,c, m_1,..., m_{t+1}$ be any non-negative integers and $F(A_1:W_1),..., F(A_t:W_t)$ be any budded ferns that satisfy all four conditions stated above. Let $p:=\sum_{i=1}^{t}{p_i}$, $n:=\sum_{i=1}^{t}{n_i}$, $w:=\sum_{i=1}^{t}{\sum_{j=1}^{r_i-1}w_j^i}$ and $m:=\sum_{i=1}^{t+1}m_i$. Suppose indices satisfy following additional conditions:\\
1) $p+n+w+m=a+b$\\
2) $p+w=n+w\leq b$. Then we have

\begin{equation}
\begin{aligned}
    &\frac{M_\odot(H_{a,b,c}^{0,\frac{b+c}{2}}(F(A_1:W_1),..., F(A_t:W_t) : m_1, m_2,..., m_t, m_{t+1}))}{M_\odot(H_{a,b,c}^{0,\frac{b+c}{2}}(B(A_1:W_1),..., B(A_t:W_t) : m_1, m_2,..., m_t, m_{t+1}))}\\
    &=\sqrt{\frac{M(H_{a,b,c}^{0,\frac{b+c}{2}}(F(A_1:W_1),..., F(A_t:W_t) : m_1, m_2,..., m_t, m_{t+1}))}{M(H_{a,b,c}^{0,\frac{b+c}{2}}(B(A_1:W_1),..., B(A_t:W_t) : m_1, m_2,..., m_t, m_{t+1}))}}\\
    &=\frac{s(X^1)}{s(Y^1)}\cdot\prod_{i=1}^{t}\Bigg[\frac{H(d_{SE}(L^i_{u_i}))H(d_{NE}(T^i))}{H(d_{SE}(T^i))H(d_{NE}(L^i_{u_i}))}\\
    &\cdot\prod_{j<u_i, j\in J_i}\frac{H(d_{SE}(L^i_j))H(d_{NE}(R^i_j))}{H(d_{SE}(R^i_j))H(d_{NE}(L^i_j))}\prod_{j\geq u_i, j\in I_i}\frac{H(d_{SE}(R^i_j))H(d_{NE}(L^i_j))}{H(d_{SE}(L^i_j))H(d_{NE}(R^i_j))}\Bigg]
\end{aligned}    
\end{equation}
\end{thm}

\section{Proof of the main results}

A region on a triangular lattice is called \textit{balanced} if it contains same number of up-pointing and down-pointing unit-triangles. Let's recall a useful result which is implicit in work of Ciucu [1] (See also Ciucu and Lai [7]).

\begin{lem}
(Region-splitting Lemma). Let $R$ be a balanced region on a triangular lattice. Assume that a subregion $S$ of $R$ satisfies the following two conditions:\\
 (1) (Seperating Condition) There is only one type of unit-triangle (either up-pointing or down-pointing) running along each side of the border between $S$ and $R-S$\\
 (2) (Balancing Condition) $S$ is balanced.
 Then
 \begin{equation}
     M(R)=M(S)M(R-S)
 \end{equation}
\end{lem}

To prove theorems in this paper, we need to simplify expressions that involves $\Delta$. For this purpose, let's recall a property of $\Delta$:

Let $X=\{x+1,x+2_,...,x+m\}$ and $Y=\{y+1,y+2,...,y+n\}$ be two sets of consecutive integers such that $x+m<y+1$. Then
\begin{equation}
\begin{aligned}
    \Delta(X,Y)=\prod_{i=1}^{m}(y-x-m+i)_n&=\prod_{i=1}^{m}\frac{(y-x+n-m+i-1)!}{(y-x-m+i-1)!}\\
    &=\frac{H(y-x-m)H(y-x+n)}{H(y-x)H(y-x+n-m)}
\end{aligned}
\end{equation}

The crucial idea of this paper is the following:\\
Each of our three main results involves the ratio of the number of tilings of two regions. For each of these two regions, we partition the set of lozenge tilings of each region according to the positions of the vertical lozenges that are bisectecd by the baseline. The partition classes obtained for the numerator and denominator are naturally paired up. Then, using Proposition 2.2. and Lemma 3.1., we verify that the ratio of the number of tilings in the corresponding partition classes does not depend on the choice of partition class (i.e. it is the same for all classes of the partition).\\

\textit{Proof of Theorem 2.1.} Let's first consider a case when $b < l \leq c$.\\
From any lozenge tiling of $H_{a,b,c}^{k,l}(X:Y)$, we will generate a pair of lozenge tiling of two trapezoidal regions with some dent on top (or bottom). If we focus on lozenges below the baseline, then the lozenges form a pentagonal region that has $b$ down-pointing unit-triangle dents on top. Among $b$ dents, $m$ of them are from the region $H_{a,b,c}^{k,l}(X:Y)$ itself, namely down-pointing unit-triangles whose bases are labeled by $y_1$, $y_2$,..., $y_m$, and remaining $(b-m)$ of them are down-pointing unit-triangles whose labels of their bases are from $[L(l)]\setminus (X\cup Y)=[a+b+k]\setminus (X\cup Y)$. Let $Z:=\{z_1, z_2,..., z_{b-m}\} \subset [L(l)]\setminus (X\cup Y)$ be a set of labels of bases of remaining $(b-m)$ dents, and let $B:=\{-|b-l|+1, -|b-l|+2,..., -1, 0\}$. Then, we can easily see that there is a natural bijection between a set of lozenge tilings of the pentagonal region having $b$ down-pointing unit-triangle dents on top and a set of lozenge tilings of a region $T(B \cup Z\cup Y)$. So from a lozenge tiling of $H_{a,b,c}^{k,l}(X:Y)$, we generate a lozenge tiling of a region $T(B \cup Z\cup Y)$.

\begin{figure}
    \centering
    \includegraphics[width=11cm]{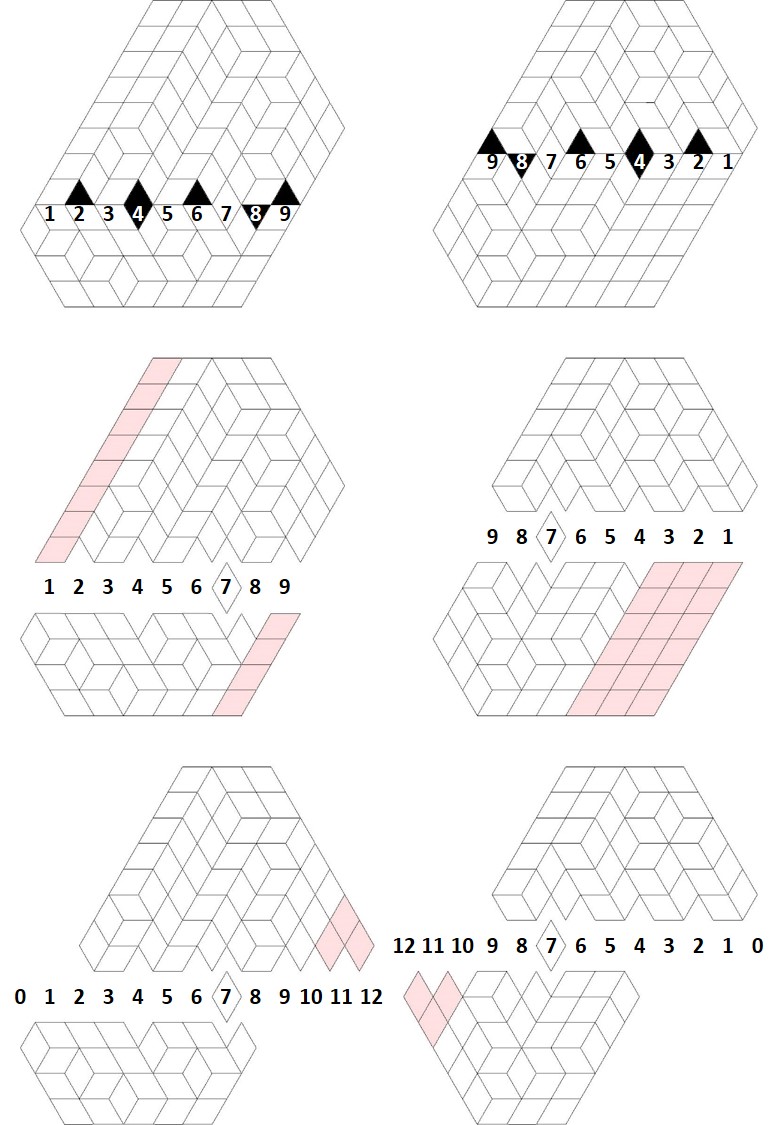}
    \caption{Correspondence between a lozenge tiling and a pair of trapezoid regions with dents (when $b<l\leq c$)}
\end{figure}

Now we return to the lozenge tiling of $H_{a,b,c}^{k,l}(X:Y)$, and focus on lozenges above the baseline. Again, they from a pentagonal region that have $(b+k)$ up-pointing unit-triangle dents on bottom. Among $(b+k)$ dents, $(m+k)$ of them are from the region $H_{a,b,c}^{k,l}(X:Y)$ itself, namely up-pointing unit-triangles whose bases are labeled by $x_1$, $x_2$,..., $x_{m+k}$ and remaining $(b-m)$ of them are up-pointing unit-triangles whose labels form a set $Z$. Let $C:=\{L(l)+1,L(l)+2,...,L(l)+|c-l|\}$. Then same observation allow us to see that there is a bijection between a set of lozenge tilings of the pentagonal region having $(b+k)$ up-pointing unit-triangle dents on bottom and a set of lozenge tilings of a region $T(Z\cup X\cup C)$. Thus, we generate a lozenge tiling of a region $T(Z\cup X\cup C)$ from a lozenge tiling of $H_{a,b,c}^{k,l}(X:Y)$.
Hence, from a lozenge tiling of $H_{a,b,c}^{k,l}(X:Y)$, we generate a pair of lozenge tiling of a region $T(B \cup Z\cup Y)$ and $T(Z\cup X\cup C)$ and this correspondence is reversible (See Figure 3.1.).  Now, we partiton a set of lozenge tiling of the region $H_{a,b,c}^{k,l}(X:Y)$ by a set $Z:=\{z_1, z_2,..., z_{b-m}\} \subset [L(l)]\setminus (X\cup Y)$ which represents labels of position of vertical lozenges on the baseline. Number of lozenge tilings of the region $H_{a,b,c}^{k,l}(X:Y)$ with $(b-m)$ vertical lozenges on the baseline whose labels of position form a set $Z:=\{z_1, z_2,..., z_{b-m}\} \subset [L(l)]\setminus (X\cup Y)$ is $H_{a,b,c}^{k,l}(Z\cup X:Z\cup Y)$. Also, by Lemma 3.1., $M(H_{a,b,c}^{k,l}(Z\cup X:Z\cup Y))$ is just a product of number of lozenge tilings of two pentagonal regions with unit-triangular dents on top (or bottom). However, numbers of lozenge tilings of two pentagonal region is same as number of lozenge tilings of regions $T(B \cup Z\cup Y)$ and $T(Z\cup X\cup C)$, respectively. Thus, we have

\begin{equation}
\begin{aligned}
    &H_{a,b,c}^{k,l}(X:Y)\\
    &=\sum_{Z=\{z_1, z_2,..., z_{b-m}\} \subseteq [L(l)]\setminus (X\cup Y)}H_{a,b,c}^{k,l}(X\cup Z:Y\cup Z)\\
    &=\sum_{Z \subseteq [L(l)]\setminus (X\cup Y)}M(T(X\cup{Z}\cup{C}))M(T(B\cup{Y}\cup{Z}))\\
    &=\sum_{Z \subseteq [L(l)]\setminus (X\cup Y)}s(X\cup{Z}\cup{C})s(B\cup{Y}\cup{Z})\\
    &=\sum_{Z\subseteq [L(l)]\setminus({{X}\cup{Y}})}{\frac{\Delta(X\cup{Z}\cup{C})}{H(b+c+k-l)}}\cdot{\frac{\Delta(B\cup{Y}\cup{Z})}{H(l)}}\\
    &=\frac{1}{H(l)\cdot H(b+c+k-l)}\cdot\sum_{Z\subseteq[L(l)]\setminus({{X}\cup{Y}})}{\Delta(X\cup{Z}\cup{C})\cdot\Delta(B\cup{Y}\cup{Z})}
\end{aligned}
\end{equation}
A lozenge tiling of ${\overline{H}_{a,b,c}^{k,b+c-l}(x_1,x_2,...,x_{m+k}:y_1,y_2,...,y_m)}$ can be also analyzed in a similar way and we can express a number of lozenge tiling of it as follows:
\begin{equation}
\begin{aligned}
    &M({\overline{H}_{a,b,c}^{k,b+c-l}(X:Y)})\\
    &=\frac{1}{H(k+l)\cdot H(b+c-l)}\cdot\sum_{Z\subseteq[L(l)]\setminus({{X}\cup{Y}})}{\Delta(B\cup{X}\cup{Z})\cdot\Delta(Y\cup{Z}\cup{C})}    
\end{aligned}
\end{equation}
where the sum is taken over all $(b-m)$ elements subset $Z\subseteq[L(l)]\setminus({{X}\cup{Y}})$.\\

However, for any $(b-m)$ elements subset $Z\subseteq[L(l)]\setminus({{X}\cup{Y}})$,
\begin{equation}
\begin{aligned}
    \frac{\Delta(X\cup{Z}\cup{C})\cdot\Delta(B\cup{Y}\cup{Z})}{\Delta(B\cup{X}\cup{Z})\cdot\Delta(Y\cup{Z}\cup{C})}&=\frac{\Delta(X)\Delta(Z)\Delta(C)\Delta(X,Z)\Delta(X,C)\Delta(Z,C)}{\Delta(B)\Delta(X)\Delta(Z)\Delta(B,X)\Delta(B,Z)\Delta(X,Z)}\\
    &~~~\cdot\frac{\Delta(B)\Delta(Y)\Delta(Z)\Delta(B,Y)\Delta(B,Z)\Delta(Y,Z)}{\Delta(Y)\Delta(Z)\Delta(C)\Delta(Y,Z)\Delta(Y,C)\Delta(Z,C)}\\
    &=\frac{\Delta(X,C)\Delta(B,Y)}{\Delta(B,X)\Delta(Y,C)}
\end{aligned}
\end{equation}
Note that this ratio does not depend on a choice of a set Z. Hence, by combining (3.3), (3,4) and (3.5), we have

\begin{equation}
\begin{aligned}
    &\frac{M({H_{a,b,c}^{k,l}(X:Y)})}{M({\overline{H}_{a,b,c}^{k,b+c-l}(X:Y)})}\\
    &=\frac{H(k+l)H(b+c-l)}{H(l)H(b+c+k-l)}\cdot\frac{\sum_{Z\subseteq[a+b+k]\setminus({{X}\cup{Y}})}{\Delta(X\cup{Z}\cup{C})\cdot\Delta(B\cup{Y}\cup{Z})}}{\sum_{Z\subseteq[a+b+k]\setminus({{X}\cup{Y}})}{\Delta(B\cup{X}\cup{Z})\cdot\Delta(Y\cup{Z}\cup{C})}}\\
    &=\frac{H(k+l)H(b+c-l)}{H(l)H(b+c+k-l)}\cdot\frac{\Delta(X,C)\Delta(B,Y)}{\Delta(B,X)\Delta(Y,C)}\\
    &=\frac{H(k+l)H(b+c-l)}{H(l)H(b+c+k-l)}\cdot\frac{\prod_{i=1}^{m+k}{(a+b+k+1-x_i)_{(c-l)}}\cdot\prod_{j=1}^{m}{(y_j)_{(l-b)}}}{\prod_{i=1}^{m+k}{(x_i)_{(l-b)}}\cdot\prod_{j=1}^{m}{(a+b+k+1-y_j)_{(c-l)}}}\\
    &=\frac{H(k+l)H(b+c-l)}{H(l)H(b+c+k-l)}\cdot\frac{\prod_{i=1}^{m+k}(x_i+l-b)_{(b-l)}(a+b+k+1-x_i)_{(c-l)}}{\prod_{j=1}^{m}(y_i+l-b)_{(b-l)}(a+b+k+1-y_j)_{(c-l)}}\\
    &=\frac{H(k+l)H(b+c-l)}{H(l)H(b+c+k-l)}\\
    &~~~\cdot\frac{\prod_{i=1}^{m+k}(x_i-b+max(b, l))_{(b-l)}(a+k+min(b, l)+1-x_i)_{(c-l)}}{\prod_{j=1}^{m}(y_i-b+max(b,l))_{(b-l)}(a+k+min(b, l)+1-y_j)_{(c-l)}}
\end{aligned}
\end{equation}

Now let's consider a case when $l \leq b$\\
Similar observation enable us to observe that $M(H_{a,b,c,k,l}(X:Y))$ can be written as a sum of $M(H_{a,b,c}^{k,l}(X\cup Z:Y\cup Z))$, where $Z=\{z_1, z_2,...,z_{l-m}\}\subset[L(l)]\setminus(X\cup Y)=[a+k+l]\setminus(X\cup Y)$ represents a set of labels of positions of vertical lozenges on the baseline. Also, by Lemma and 3.1. and same argument as we used in previous case, $M(H_{a,b,c}^{k,l}(X\cup Z:Y\cup Z))$ is equal to a product of $M(T(B\cup{X}\cup{Z}\cup{C}))$ and $M(T(Y\cup{Z}))$, where $B=\{-|b-l|+1, -|b-l|+2,..., -1, 0\}$ and $C=\{L(l)+1, L(l)+2,..., L(l)+|c-l|\}$. Hence

\begin{figure}
    \centering
    \includegraphics[width=11cm]{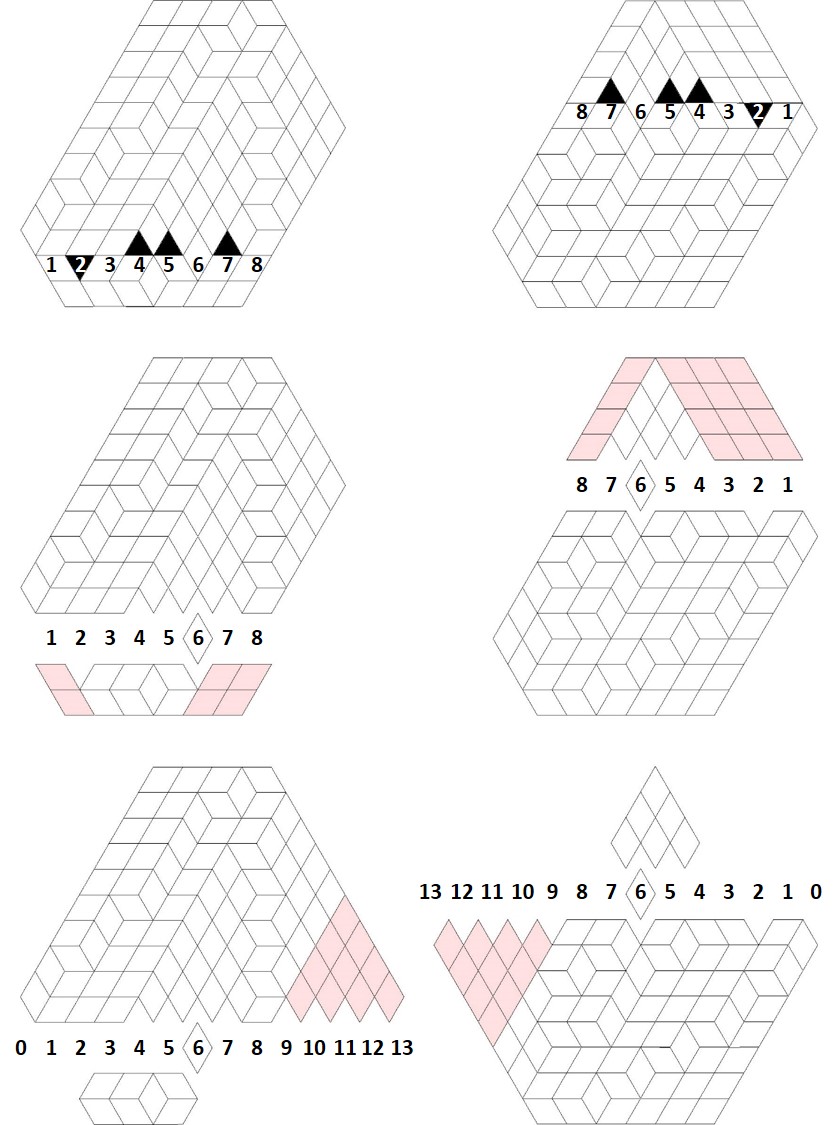}
    \caption{Correspondence between a lozenge tiling and a pair of trapezoid regions with dents (when $l\leq b$)}
\end{figure}

\begin{equation}
H_{a,b,c}^{k,l}(Z\cup X:Z\cup Y)=\frac{\Delta(B\cup{X}\cup{Z}\cup{C})}{H(b+c+k-l)}\cdot\frac{\Delta(Y\cup{Z})}{H(l)}
\end{equation}
If we sum over every $(l-m)$ element set $Z\subseteq[L(l)]\setminus({{X}\cup{Y}})$, then we have a representation of number of lozenge tiling of $H_{a,b,c}^{k,l}(X:Y)$ as follows:

\begin{equation}
\begin{aligned}
    &M(H_{a,b,c}^{k,l}(X:Y))\\
    &=\frac{1}{H(l)\cdot H(b+c+k-l)}\cdot\sum_{Z\subseteq[L(l)]\setminus({{X}\cup{Y}})}{\Delta(B\cup{X}\cup{Z}\cup{C})\cdot\Delta(Y\cup{Z})}
\end{aligned}
\end{equation}
By same observation, we can represent a number of lozenge tiling of a hexagon $\overline{H}_{a,b,c}^{k,b+c-l}(X:Y)$ as follows:
\begin{equation}
\begin{aligned}
    &M(\overline{H}_{a,b,c}^{k,b+c-l}(X:Y))\\
    &=\frac{1}{H(k+l)\cdot H(b+c-l)}\cdot\sum_{Z\subseteq[L(l)]\setminus({{X}\cup{Y}})}{\Delta(X\cup{Z})\cdot\Delta(B\cup{Y}\cup{Z}\cup{C})}
\end{aligned}
\end{equation}
Now, we observe a ratio $\frac{\Delta(B\cup{X}\cup{Z}\cup{C})\cdot\Delta(Y\cup{Z})}{\Delta(X\cup{Z})\cdot\Delta(B\cup{Y}\cup{Z}\cup{C})}$ for any subset $Z\subseteq[L(l)]\setminus({{X}\cup{Y}})$ with $(l-m)$ elements:
\begin{equation}
\begin{aligned}
    &\frac{\Delta(B\cup{X}\cup{Z}\cup{C})\cdot\Delta(Y\cup{Z})}{\Delta(X\cup{Z})\cdot\Delta(B\cup{Y}\cup{Z}\cup{C})}\\
    &=\frac{\Delta(B)\Delta(X)\Delta(Z)\Delta(C)\Delta(B,X)\Delta(B,Z)\Delta(B,C)\Delta(X,Z)\Delta(X,C)\Delta(Z,C)}{\Delta(X)\Delta(Z)\Delta(X,Z)}\\
    &\cdot\frac{\Delta(Y)\Delta(Z)\Delta(Y,Z)}{\Delta(B)\Delta(Y)\Delta(Z)\Delta(C)\Delta(B,Y)\Delta(B,Z)\Delta(B,C)\Delta(Y,Z)\Delta(Y,C)\Delta(Z,C)}\\
    &=\frac{\Delta(B,X)\Delta(X,C)}{\Delta(B,Y)\Delta(Y,C)}\\
    &=\frac{\prod_{i=1}^{m+k}{(x_i)_{(b-l)}\cdot(a+b+k+1-x_i)}_{(c-l)}}{\prod_{j=1}^{m}{(y_j)_{(b-l)}\cdot(a+b+k+1-y_j)_{(c-l)}}}
\end{aligned}
\end{equation}

Note that this ratio does not depend on a choice of a set Z. Hence, by combining (3.8), (3.9) and (3.10), we have

\begin{equation}
\begin{aligned}
    &\frac{M({H_{a,b,c}^{k,l}(X:Y)})}{M({\overline{H}_{a,b,c}^{k,b+c-l}(X:Y)})}\\
    &=\frac{H(k+l)H(b+c-l)}{H(l)H(b+c+k-l)}\cdot\frac{\sum_{Z\subseteq[a+k+l]\setminus({{X}\cup{Y}})}{\Delta(B\cup{X}\cup{Z}\cup{C})\cdot\Delta(Y\cup{Z})}}{\sum_{Z\subseteq[a+k+l]\setminus({{X}\cup{Y}})}{\Delta(X\cup{Z})\cdot\Delta(B\cup{Y}\cup{Z}\cup{C})}}\\
    &=\frac{H(k+l)H(b+c-l)}{H(l)H(b+c+k-l)}\cdot\frac{\prod_{i=1}^{m+k}{(x_i)_{(b-l)}\cdot(a+b+k+1-x_i)}_{(c-l)}}{\prod_{j=1}^{m}{(y_j)_{(b-l)}\cdot(a+b+k+1-y_j)_{(c-l)}}}\\
    &=\frac{H(k+l)H(b+c-l)}{H(l)H(b+c+k-l)}\\
    &~~~\cdot\frac{\prod_{i=1}^{m+k}{(x_i-b+max(b, l))_{(b-l)}\cdot(a+k+max(b, l)+1-x_i)}_{(c-l)}}{\prod_{j=1}^{m}{(y_j-b+max(b, l))_{(b-l)}\cdot(a+k+max(b, l)+1-y_j)_{(c-l)}}}
\end{aligned}
\end{equation}

The case when $c < l \leq b+c$ can be proved similarly as we did for the case when $l \leq b$. Hence, the theorem has been proved. $\square$ \\

\begin{figure}
    \centering
    \includegraphics[width=1\textwidth]{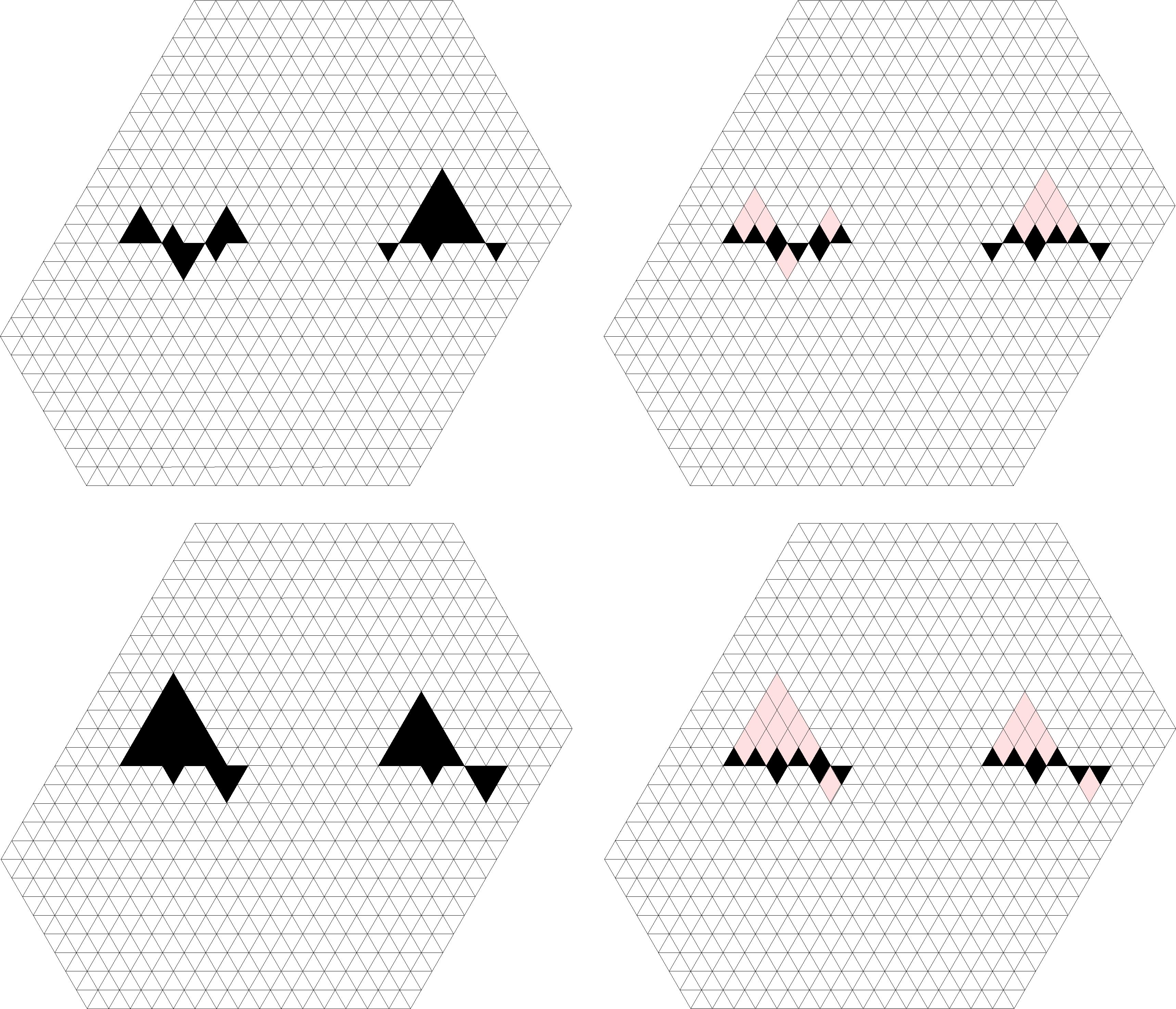}
    \caption{}
\end{figure}

\textit{Proof of Theorem 2.3}. Again, let's consider a case when $b < l \leq c$ first. If we compare two regions $H_{a,b,c}^{k,l}(F(A_1:W_1),..., F(A_t:W_t) : m_1, ..., m_{t+1})$ and $H_{a,b,c}^{k,l}(F_{br}(A_1:W_1),..., F_{br}(A_t:W_t) : m_1,..., m_{t+1})$, two regions are different by equilateral triangles with zig-zag horizontal boundary. However, in any lozenge tiling of a region $H_{a,b,c}^{k,l}(F_{br}(A_1:W_1),..., F_{br}(A_t:W_t) : m_1, m_2,..., m_t, m_{t+1})$, those regions are forced to be tiled by vertical lozenges (See Figure 3.3). Hence, two regions have same number of lozenge tiling. Similarly, two regions $H_{a,b,c}^{k,l}(B(A_1:W_1),..., B(A_t:W_t) : m_1,..., m_{t+1})$ and $H_{a,b,c}^{k,l}(B_{br}(A_1:W_1),..., B_{br}(A_t:W_t) : m_1,..., m_{t+1})$ have same number of lozenge tilings. Thus we have
\begin{equation}
\begin{aligned}
    &\frac{M(H_{a,b,c}^{k,l}(F(A_1:W_1),..., F(A_t:W_t) : m_1,..., m_{t+1}))}{M(H_{a,b,c}^{k,l}(B(A_1:W_1),..., B(A_t:W_t) : m_1,..., m_{t+1}))}\\
    &=\frac{M(H_{a,b,c}^{k,l}(F_{br}(A_1:W_1),..., F_{br}(A_t:W_t) : m_1,..., m_{t+1}))}{M(H_{a,b,c}^{k,l}(B_{br}(A_1:W_1),..., B_{br}(A_t:W_t) : m_1,..., m_{t+1}))}
\end{aligned}
\end{equation}

For $i\in[t]$, $j\in [r_i]$, let $X^i_j=\{d_{NW}(L^i_j)+1, d_{NW}(L^i_j)+2,..., d_{NW}(R^i_j) (=d_{NW}(L^i_j)+a^i_j)\}$, $V_i=\{d_{NW}(L^i_j)+1, d_{NW}(L^i_j)+2,..., d_{NW}(T^i) (=d_{NW}(L^i_j)+v_i)\}$ and $\overline{V_i}=X^i_{u_i}\setminus V_i=\{d_{NW}(T^i)+1,..., d_{NW}(R^i_j)\}$. Then $X^1=\cup_{i=1}^{t}\cup_{j\in I_i}X^i_j$, $X^2=\cup_{i=1}^{t}\cup_{j\in J_i}X^i_j$, $Y^1=\cup_{i=1}^{t}((\cup_{j=1}^{u_i-1}X^i_j)\cup V_i)$ and $Y^2=\cup_{i=1}^{t}(\overline{V_i}\cup (\cup_{j=u_i+1}^{r_i}X^i_j))$.

By same observation as we did in the Theorem 1, Lozenge tiling of a hexagonal region $H_{a,b,c}^{k,l}(F_{br}(A_1:W_1),..., F_{br}(A_t:W_t) : m_1, m_2,..., m_t, m_{t+1})$ can be partitioned according to $(b-n-w)$ vertical unit-lozenges that are bisected by the $l$-th horizontal line. Let $H_{a,b,c}^{k,l}(F_{br}(A_1:W_1),..., F_{br}(A_t:W_t) : m_1, m_2,..., m_t, m_{t+1} : z_1,..., z_{b-n-w})$ be a region obtained from $H_{a,b,c}^{k,l}(F_{br}(A_1:W_1),..., F_{br}(A_t:W_t) : m_1, m_2,..., m_t, m_{t+1})$ by removing $(b-n-w)$ unit-lozenges that are bisected by segments on $l$-th horizontal line whose labels are elements of a set $Z=\{z_1,z_2,...,z_{b-n-w}\}$. Then, by same argument as we used in proof of Theorem 2.1., we have
\begin{equation}
\begin{aligned}
    &M(H_{a,b,c}^{k,l}(F_{br}(A_1:W_1),..., F_{br}(A_t:W_t) : m_1,..., m_{t+1} : z_1,..., z_{b-n-w}))\\
    &=s(Z\cup{X^1}\cup W\cup C)\cdot s(B\cup Z\cup{X^2}\cup W)\\
    &=\frac{\Delta(Z\cup{X^1}\cup W\cup C)}{H(b+c+k-l)}\cdot\frac{\Delta(B\cup Z\cup{X^2}\cup W)}{H(l)}
\end{aligned}
\end{equation}
where $B=\{-|b-l|+1,..., -1, 0\}$ and $C=\{L(l)+1, L(l)+2,..., L(l)+|c-l|\}$.
If we sum over every $(b-n-w)$ element set $Z\subset[L(l)]\setminus(X^1\cup X^2\cup W)$, then we have a representation of number of lozenge tilings of $H_{a,b,c}^{k,l}(F_{br}(A_1:W_1),..., F_{br}(A_t:W_t) : m_1,..., m_{t+1})$ as follows:
\begin{equation}
\begin{aligned}
    &M(H_{a,b,c}^{k,l}(F_{br}(A_1:W_1),..., F_{br}(A_t:W_t) : m_1,..., m_{t+1}))\\
    &=\sum_{Z}\frac{\Delta(Z\cup{X^1}\cup W\cup C)\cdot\Delta(B\cup Z\cup{X^2}\cup W)}{H(l)\cdot H(b+c+k-l)}\\
    &=\frac{\sum_{Z}\Delta(Z\cup{X^1}\cup W\cup C)\cdot\Delta(B\cup Z\cup{X^2}\cup W)}{H(l)\cdot H(b+c+k-l)}\\
\end{aligned}
\end{equation}

Similarly, a number of lozenge tilings of
$H_{a,b,c}^{k,l}(B_{br}(A_1:W_1),..., B_{br}(A_t:W_t) : m_1,..., m_{t+1})$ can be expressed as follow:
\begin{equation}
\begin{aligned}
    &M(H_{a,b,c}^{k,l}(B_{br}(A_1:W_1),..., B_{br}(A_t:W_t) : m_1,..., m_{t+1}))\\
    &=\sum_{Z}\frac{\Delta(Z\cup{Y^1}\cup W\cup C)\cdot\Delta(B\cup Z\cup{Y^2}\cup W)}{H(l)\cdot H(b+c+k-l)}\\
    &=\frac{\sum_{Z}\Delta(Z\cup{Y^1}\cup W\cup C)\cdot\Delta(B\cup Z\cup{Y^2}\cup W)}{H(l)\cdot H(b+c+k-l)}\\
\end{aligned}
\end{equation}

Now, let's observe a ratio $\frac{\Delta(Z\cup{X^1}\cup W\cup C)\cdot\Delta(B\cup Z\cup{X^2}\cup W)}{\Delta(Z\cup{Y^1}\cup W\cup C)\cdot\Delta(B\cup Z\cup{Y^2}\cup W)}$ for any set $Z\subset[L(l)]\setminus(X^1\cup X^2\cup W)$ with ($b-n-w$) elements:
\begin{equation}
\begin{aligned}
    &\frac{\Delta(Z\cup{X^1}\cup W\cup C)\cdot\Delta(B\cup Z\cup{X^2}\cup W)}{\Delta(Z\cup{Y^1}\cup W\cup C)\cdot\Delta(B\cup Z\cup{Y^2}\cup W)}\\
    &=\frac{\Delta(Z)\Delta(X^1)\Delta(W)\Delta(C)\Delta(Z,X^1)\Delta(Z,W)\Delta(Z,C)\Delta(X^1,W)\Delta(X^1,C)\Delta(W,C)}{\Delta(Z)\Delta(Y^1)\Delta(W)\Delta(C)\Delta(Z,Y^1)\Delta(Z,W)\Delta(Z,C)\Delta(Y^1,W)\Delta(Y^1,C)\Delta(W,C)}\\
    &\cdot\frac{\Delta(B)\Delta(Z)\Delta(X^2)\Delta(W)\Delta(B,Z)\Delta(B,X^2)\Delta(B,W)\Delta(Z,X^2)\Delta(Z,W)\Delta(X^2,W)}{\Delta(B)\Delta(Z)\Delta(Y^2)\Delta(W)\Delta(B,Z)\Delta(B,Y^2)\Delta(B,W)\Delta(Z,Y^2)\Delta(Z,W)\Delta(Y^2,W)}\\
    &=\frac{\Delta(X^1)\Delta(X^2)\Delta(X^1,C)\Delta(B,X^2)}{\Delta(Y^1)\Delta(Y^2)\Delta(Y^1,C)\Delta(B,Y^2)}\\
    &=\frac{s(X^1)s(X^2)}{s(Y^1)s(Y^2)}\cdot\frac{\Delta(X^1,C)}{\Delta(Y^1,C)}\cdot\frac{\Delta(B,X^2)}{\Delta(B,Y^2)}
\end{aligned}
\end{equation}
In above simplification, we use a fact that $X^1 \cup X^2=Y^1 \cup Y^2$, which implies $\Delta(Z,X^1)\Delta(Z,X^2)=\Delta(Z,Y^1)\Delta(Z,Y^2)$ and $\Delta(X^1,W)\Delta(X^2,W)=\Delta(Y^1,W)\Delta(Y^2,W)$. Note that what we get does not depend on a choice of a set $Z$. Hence, by (3.14), (3.15), (3.16) and (3.17), we have
\begin{equation}
\begin{aligned}
    &\frac{M(H_{a,b,c}^{k,l}(F(A_1:W_1),..., F(A_t:W_t) : m_1, m_2,..., m_t, m_{t+1}))}{M(H_{a,b,c}^{k,l}(B(A_1:W_1),..., B(A_t:W_t) : m_1, m_2,..., m_t, m_{t+1}))}\\
    &=\frac{s(X^1)s(X^2)}{s(Y^1)s(Y^2)}\cdot\frac{\Delta(X^1,C)}{\Delta(Y^1,C)}\cdot\frac{\Delta(B,X^2)}{\Delta(B,Y^2)}
\end{aligned}
\end{equation}

Since $X^1=\cup_{i=1}^{t}\cup_{j\in I_i}X^i_j$ and $Y^1=\cup_{i=1}^{t}((\cup_{j=1}^{u_i-1}X^i_j)\cup V_i)$,
\begin{equation}
\begin{aligned}
    \frac{\Delta(X^1,C)}{\Delta(Y^1,C)}&=\frac{\prod_{i=1}^{t}\prod_{j\in I_i}\Delta(X^i_j, C)}{\prod_{i=1}^{t}((\prod_{j=1}^{u_i-1}\Delta(X^i_j, C))\cdot\Delta(V_i, C))}\\
    &=\prod_{i=1}^{t}\Bigg[\frac{1}{\Delta(V_i, C)}\prod_{j<u_i, j\in J_i}\frac{1}{\Delta(X^i_j, C)}\prod_{j\geq u_i, j\in I_i}\Delta(X^i_j, C)\Bigg]
\end{aligned}
\end{equation}

However, by (3.2), we have
\begin{equation}
\begin{aligned}
    \Delta(V_i, C)&=\frac{H(L(l)-d_{NW}(L^i_{u_i})-v_i)H(L(l)-d_{NW}(L^i_{u_i})+|c-l|)}{H(L(l)-d_{NW}(L^i_{u_i}))H(L(l)-d_{NW}(L^i_{u_i})+|c-l|-v_i)}\\
    &=\frac{H(d_{SE}(T^i))H(d_{NE}(L^i_{u_i}))}{H(d_{SE}(L^i_{u_i}))H(d_{NE}(T^i))}
\end{aligned}
\end{equation}
and
\begin{equation}
\begin{aligned}
    \Delta(X^i_j, C)&=\frac{H(L(l)-d_{NW}(L^i_{u_i})-a^i_j)H(L(l)-d_{NW}(L^i_{u_i})+|c-l|)}{H(L(l)-d_{NW}(L^i_{u_i}))H(L(l)-d_{NW}(L^i_{u_i})+|c-l|-a^i_j)}\\
    &=\frac{H(d_{SE}(R^i_j))H(d_{NE}(L^i_j))}{H(d_{SE}(L^i_j))H(d_{NE}(R^i_j))}
\end{aligned}
\end{equation}
Hence, by (3.18), (3.19) and (3.20), we have
\begin{equation}
\begin{aligned}
    \frac{\Delta(X^1,C)}{\Delta(Y^1,C)}=&\prod_{i=1}^{t}\Bigg[\frac{H(d_{SE}(L^i_{u_i}))H(d_{NE}(T^i))}{H(d_{SE}(T^i))H(d_{NE}(L^i_{u_i}))}\\
    &\cdot\prod_{j<u_i, j\in J_i}\frac{H(d_{SE}(L^i_j))H(d_{NE}(R^i_j))}{H(d_{SE}(R^i_j))H(d_{NE}(L^i_j))}\prod_{j\geq u_i, j\in I_i}\frac{H(d_{SE}(R^i_j))H(d_{NE}(L^i_j))}{H(d_{SE}(L^i_j))H(d_{NE}(R^i_j))}\Bigg]
\end{aligned}
\end{equation}

Also, since $X^2=\cup_{i=1}^{t}\cup_{j\in J_i}X^i_j$ and $Y^2=\cup_{i=1}^{t}(\overline{V_i} \cup (\cup_{j=u_i+1}^{r_i}X^i_j))$,
\begin{equation}
\begin{aligned}
    \frac{\Delta(B,X^2)}{\Delta(B,Y^2)}&=\frac{\prod_{i=1}^{t}\prod_{j\in J_i}\Delta(B, X^i_j)}{\prod_{i=1}^{t}(\Delta(B, \overline{V_i})\cdot\prod_{j=u_i+1}^{r_i}\Delta(B, X^i_j))}\\
    &=\prod_{i=1}^{t}\Bigg[\frac{\Delta(B, X^i_{u_i})}{\Delta(B, \overline{V_i})}\prod_{j<u_i, j\in J_i}\Delta(B, X^i_j)\prod_{j\geq u_i, j\in I_i}\frac{1}{\Delta(B, X^i_j)}\Bigg]\\
    &=\prod_{i=1}^{t}\Bigg[\Delta(B, V_i)\prod_{j<u_i, j\in J_i}\Delta(B, X^i_j)\prod_{j\geq u_i, j\in I_i}\frac{1}{\Delta(B, X^i_j)}\Bigg]
\end{aligned}
\end{equation}

Again, by (3.2), we have
\begin{equation}
\begin{aligned}
    \Delta(B, V_i)&=\frac{H(d_{NW}(L^i_{u_i}))H(d_{NW}(L^i_{u_i})+|l-b|+v_i)}{H(d_{NW}(L^i_{u_i})+|l-b|)H(d_{NW}(L^i_{u_i})+v_i)}\\
    &=\frac{H(d_{NW}(L^i_{u_i}))H(d_{SW}(T^i))}{H(d_{SW}(L^i_{u_i}))H(d_{NW}(T^i))}\\
    &=\frac{H(d_{NW}(L^i_{u_i}))H(d_{SW}(T^i))}{H(d_{NW}(T^i))H(d_{SW}(L^i_{u_i}))}
\end{aligned}
\end{equation}
and
\begin{equation}
\begin{aligned}
    \Delta(B, X^i_j)&=\frac{H(d_{NW}(L^i_{u_i}))H(d_{NW}(L^i_{u_i})+|l-b|+a^i_j)}{H(d_{NW}(L^i_{u_i})+|l-b|)H(d_{NW}(L^i_{u_i})+a^i_j)}\\
    &=\frac{H(d_{NW}(L^i_j))H(d_{SW}(R^i_j))}{H(d_{SW}(L^i_j))H(d_{NW}(R^i_j))}\\
    &=\frac{H(d_{NW}(L^i_j))H(d_{SW}(R^i_j))}{H(d_{NW}(R^i_j))H(d_{SW}(L^i_j))}
\end{aligned}
\end{equation}

Hence, by (3.22), (3.23) and (3.24), we have
\begin{equation}
\begin{aligned}
    \frac{\Delta(B,X^2)}{\Delta(B,Y^2)}=&\prod_{i=1}^{t}\Bigg[\frac{H(d_{NW}(L^i_{u_i}))H(d_{SW}(T^i))}{H(d_{NW}(T^i))H(d_{SW}(L^i_{u_i}))}\\
    &\cdot\prod_{j<u_i, j\in J_i}\frac{H(d_{NW}(L^i_j))H(d_{SW}(R^i_j))}{H(d_{NW}(R^i_j))H(d_{SW}(L^i_j))}\prod_{j\geq u_i, j\in I_i}\frac{H(d_{NW}(R^i_j))H(d_{SW}(L^i_j))}{H(d_{NW}(L^i_j))H(d_{SW}(R^i_j))}\Bigg]
\end{aligned}
\end{equation}

Thus, by (3.17), (3.21) and (3.25),
\begin{equation}
\begin{aligned}
    &\frac{M(H_{a,b,c}^{k,l}(F(A_1:W_1),..., F(A_t:W_t) : m_1, m_2,..., m_t, m_{t+1}))}{M(H_{a,b,c}^{k,l}(B(A_1:W_1),..., B(A_t:W_t) : m_1, m_2,..., m_t, m_{t+1}))}\\
    &=\frac{s(X^1)s(X^2)}{s(Y^1)s(Y^2)}\\
    &\cdot\prod_{i=1}^{t}\Bigg[\frac{H(d_{SE}(L^i_{u_i}))H(d_{NE}(T^i))H(d_{NW}(L^i_{u_i}))H(d_{SW}(T^i))}{H(d_{SE}(T^i))H(d_{NE}(L^i_{u_i}))H(d_{NW}(T^i))H(d_{SW}(L^i_{u_i}))}\\
    &\cdot \prod_{j<u_i, j\in J_i}\frac{H(d_{SE}(L^i_j))H(d_{NE}(R^i_j))H(d_{NW}(L^i_j))H(d_{SW}(R^i_j))}{H(d_{SE}(R^i_j))H(d_{NE}(L^i_j))H(d_{NW}(R^i_j))H(d_{SW}(L^i_j))}\\
    &\cdot \prod_{j\geq u_i, j\in I_i}\frac{H(d_{SE}(R^i_j))H(d_{NE}(L^i_j))H(d_{NW}(R^i_j))H(d_{SW}(L^i_j))}{H(d_{SE}(L^i_j))H(d_{NE}(R^i_j))H(d_{NW}(L^i_j))H(d_{SW}(R^i_j))}\Bigg]
\end{aligned}
\end{equation}

Now, let's consider a case when $l \leq b$.\\
For $i\in[t]$, $j\in [r_i]$, let $X^i_j=\{d_{SW}(L^i_j)+1, d_{SW}(L^i_j)+2,..., d_{SW}(R^i_j) (=d_{SW}(L^i_j)+a^i_j)\}$, $V_i=\{d_{SW}(L^i_j)+1, d_{SW}(L^i_j)+2,..., d_{SW}(T^i) (=d_{SW}(L^i_j)+v_i)\}$ and $\overline{V_i}=X^i_{u_i}\setminus V_i=\{d_{SW}(T^i)+1,..., d_{SW}(R^i_j)\}$. Then $X^1=\cup_{i=1}^{t}\cup_{j\in I_i}X^i_j$, $X^2=\cup_{i=1}^{t}\cup_{j\in J_i}X^i_j$, $Y^1=\cup_{i=1}^{t}((\cup_{j=1}^{u_i-1}X^i_j)\cup V_i)$ and $Y^2=\cup_{i=1}^{t}(\overline{V_i}\cup (\cup_{j=u_i+1}^{r_i}X^i_j))$.\\
By same argument, the ratio can be expressed as follows:
\begin{equation}
\begin{aligned}
    &\frac{M(H_{a,b,c}^{k,l}(F(A_1:W_1),..., F(A_t:W_t) : m_1, m_2,..., m_t, m_{t+1}))}{M(H_{a,b,c}^{k,l}(B(A_1:W_1),..., B(A_t:W_t) : m_1, m_2,..., m_t, m_{t+1}))}\\
    &=\frac{s(X^1)s(X^2)}{s(Y^1)s(Y^2)}\cdot\frac{\Delta(X^1,C)}{\Delta(Y^1,C)}\cdot\frac{\Delta(B,X^1)}{\Delta(B,Y^1)}
\end{aligned}
\end{equation}
where $B=\{-|b-l|+1,..., -1, 0\}$ and $C=\{L(l)+1, L(l)+2,..., L(l)+|c-l|\}$.
However, we know that
\begin{equation}
\begin{aligned}
    \frac{\Delta(X^1,C)}{\Delta(Y^1,C)}&=\frac{\prod_{i=1}^{t}\prod_{j\in I_i}\Delta(X^i_j, C)}{\prod_{i=1}^{t}((\prod_{j=1}^{u_i-1}\Delta(X^i_j, C))\cdot\Delta(V_i, C))}\\
    &=\prod_{i=1}^{t}\Bigg[\frac{1}{\Delta(V_i, C)}\prod_{j<u_i, j\in J_i}\frac{1}{\Delta(X^i_j, C)}\prod_{j\geq u_i, j\in I_i}\Delta(X^i_j, C)\Bigg]\\
\end{aligned}
\end{equation}
Also, we have
\begin{equation}
\begin{aligned}
    \frac{\Delta(B,X^1)}{\Delta(B,Y^1)}&=\frac{\prod_{i=1}^{t}\prod_{j\in I_i}\Delta(B, X^i_j)}{\prod_{i=1}^{t}((\prod_{j=1}^{u_i-1}\Delta(B, X^i_j))\cdot\Delta(B, V_i))}\\
    &=\prod_{i=1}^{t}\Bigg[\frac{1}{\Delta(B, V_i)}\prod_{j<u_i, j\in J_i}\Delta(B, X^i_j)\prod_{j\geq u_i, j\in I_i}\frac{1}{\Delta(B, X^i_j)}\Bigg]\\
\end{aligned}
\end{equation}

However, by (3.2), we have

\begin{equation}
\begin{aligned}
    \Delta(B, X^i_j)&=\frac{H(d_{SW}(L^i_j))H(d_{SW}(L^i_j)+|b-l|+a^i_j)}{H(d_{SW}(L^i_j)+|b-l|)H(d_{SW}(L^i_j)+a^i_j)}\\
    &=\frac{H(d_{SW}(L^i_j))H(d_{NW}(R^i_j))}{H(d_{SW}(R^i_j))H(d_{NW}(L^i_j))}\\
\end{aligned}
\end{equation}

\begin{equation}
\begin{aligned}
    \Delta(X^i_j, C)&=\frac{H(L(l)-d_{SW}(L^i_j)-a^i_j)H(L(l)-d_{SW}(L^i_j)+|c-l|)}{H(L(l)-d_{SW}(L^i_j))H(L(l)-d_{SW}(L^i_j)+|c-l|-a^i_j)}\\
    &=\frac{H(d_{SE}(R^i_j))H(d_{NE}(L^i_j))}{H(d_{SE}(L^i_j))H(d_{NE}(R^i_j))}\\
\end{aligned}
\end{equation}

and similarly
\begin{equation}
\begin{aligned}
    \Delta(B, V_i)&=\frac{H(d_{SW}(L^i_{u_i}))H(d_{NW}(T^i))}{H(d_{SW}(T^i))H(d_{NW}(L^i_{u_i}))}, \Delta(V_i, C)=\frac{H(d_{SE}(T_i))H(d_{NE}(L^i_{u_i}))}{H(d_{SE}(L^i_{u_i}))H(d_{NE}(T_i))}
\end{aligned}
\end{equation}

Thus, by (3.27)-(3.32), we have
\begin{equation}
\begin{aligned}
    &\frac{M(H_{a,b,c}^{k,l}(F(A_1:W_1),..., F(A_t:W_t) : m_1, m_2,..., m_t, m_{t+1}))}{M(H_{a,b,c}^{k,l}(B(A_1:W_1),..., B(A_t:W_t) : m_1, m_2,..., m_t, m_{t+1}))}\\
    &=\frac{s(X^1)s(X^2)}{s(Y^1)s(Y^2)}\\
    &\cdot\prod_{i=1}^{t}\Bigg[\frac{H(d_{SE}(L^i_{u_i}))H(d_{NE}(T^i))H(d_{NW}(L^i_{u_i}))H(d_{SW}(T^i))}{H(d_{SE}(T^i))H(d_{NE}(L^i_{u_i}))H(d_{NW}(T^i))H(d_{SW}(L^i_{u_i}))}\\
    &\cdot \prod_{j<u_i, j\in J_i}\frac{H(d_{SE}(L^i_j))H(d_{NE}(R^i_j))H(d_{NW}(L^i_j))H(d_{SW}(R^i_j))}{H(d_{SE}(R^i_j))H(d_{NE}(L^i_j))H(d_{NW}(R^i_j))H(d_{SW}(L^i_j))}\\
    &\cdot \prod_{j\geq u_i, j\in I_i}\frac{H(d_{SE}(R^i_j))H(d_{NE}(L^i_j))H(d_{NW}(R^i_j))H(d_{SW}(L^i_j))}{H(d_{SE}(L^i_j))H(d_{NE}(R^i_j))H(d_{NW}(L^i_j))H(d_{SW}(R^i_j))}\Bigg]
\end{aligned}
\end{equation}

The case when $c <l$ can be proved similarly as we did for the case when $l \leq b$. Hence, the Theorem 2 has been proved. $\square$\\

\begin{figure}
    \centering
    \includegraphics[width=1\textwidth]{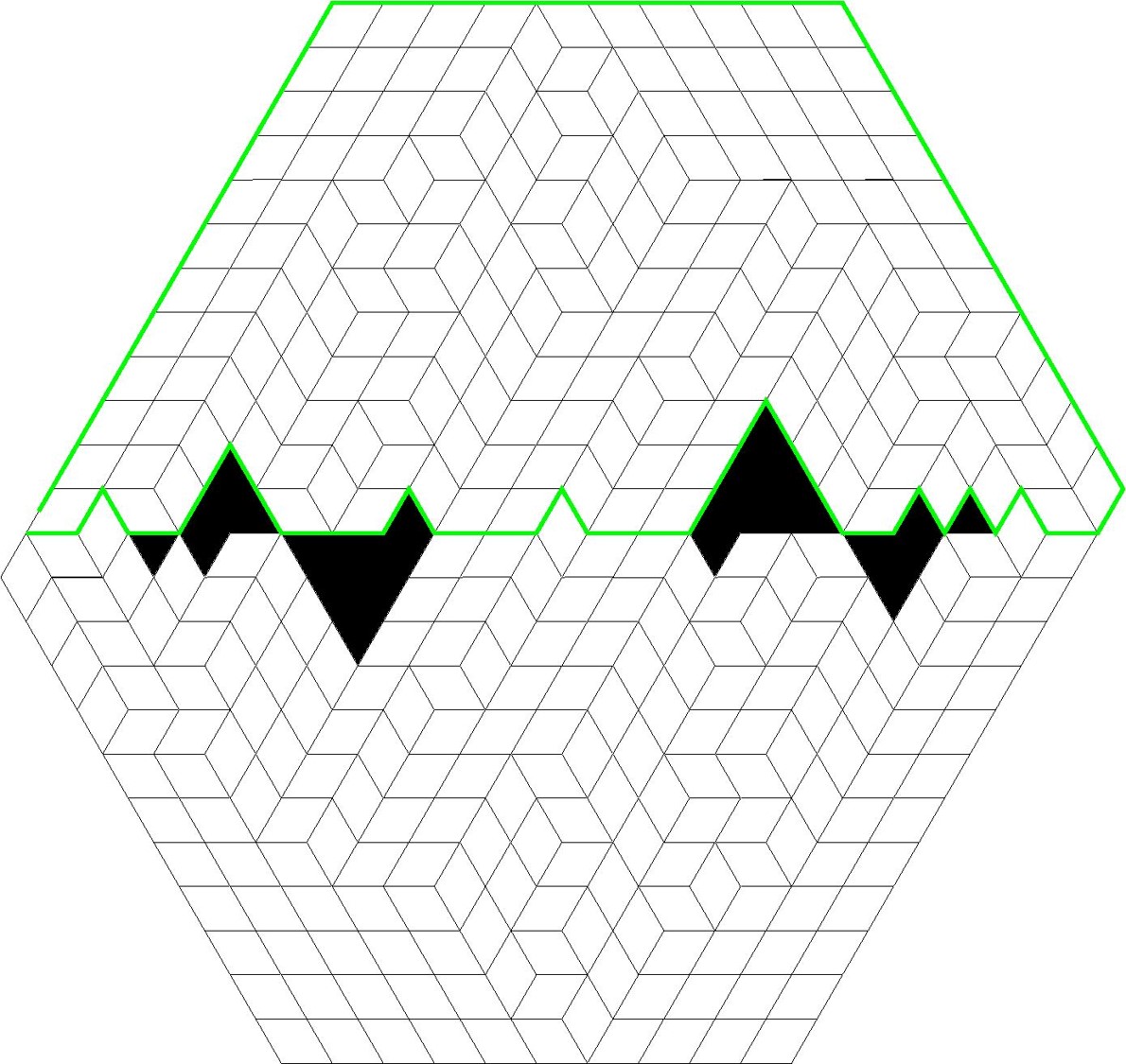}
    \caption{Cyclically symmetric lozenge tiling of a region $H_{10,11,13}^{0,12}(F(-1,1,-2,0:1,0,1), F(0,2,-1,1:1,0,1) : 2, 5, 2)$}
\end{figure}

\textit{Proof of Theorem 2.4}.
Let's use the same notation as we used in the proof of Theorem 2.3. Like proof of previous theorems, we label the baseline by $1, 2,..., L(\frac{b+c}{2})$ from left to right.
Note that in this case, sets $X^1$, $X^2$, $Y^1$, $Y^2$ and $W$ satisfy $X^2=\{L(\frac{b+c}{2})+1-x|x\in X^1\}$, $Y^2=\{L(\frac{b+c}{2})+1-y|y\in Y^1\}$ and $W=\{L(\frac{b+c}{2})+1-w|w\in W\}$ because the region is centrally symmetric. Crucial observation is that \textbf{a centrally symmetric lozenge tiling of the region is uniquely determined by lozenges below (or above) the horizontal line (See Figure 3.4)} .

Hence, by combining this observation and same argument that we have used in the proof of previous theorems, we have

\begin{equation}
\begin{aligned}
    &M_\odot(H_{a,b,c}^{0,\frac{b+c}{2}}(F(A_1:W_1),..., F(A_t:W_t) : m_1, m_2,..., m_t, m_{t+1}))\\
    &=\sum_{Z}\frac{\Delta(Z\cup X^1 \cup W \cup C)}{H(\frac{b+c}{2})}\\
    &=\frac{\sum_{Z}\Delta(Z\cup X^1\cup W \cup C)}{H(\frac{b+c}{2})}
\end{aligned}    
\end{equation}
where the sum is taken over all sets $Z\subset[L(\frac{b+c}{2})]\setminus(X^1\cup X^2\cup W)$ with $(b-n-w)$ elements that satisfies $Z=\{L(\frac{b+c}{2})+1-z|z\in Z\}$

Similarly, number of centrally symmetric lozenge tiling of the region $H_{a,b,c}^{0,\frac{b+c}{2}}(B(A_1:W_1),..., B(A_t:W_t) : m_1, m_2,..., m_t, m_{t+1})$ can be written as follows:
\begin{equation}
\begin{aligned}
    &M_\odot(H_{a,b,c}^{0,\frac{b+c}{2}}(B(A_1:W_1),..., B(A_t:W_t) : m_1, m_2,..., m_t, m_{t+1}))\\
    &=\sum_{Z}\frac{\Delta(Z\cup Y^1\cup W \cup C)}{H(\frac{b+c}{2})}\\
    &=\frac{\sum_{Z}\Delta(Z\cup Y^1\cup W \cup C)}{H(\frac{b+c}{2})}
\end{aligned}
\end{equation}

Again, the sum is taken over all sets $Z\subset[L(\frac{b+c}{2})]\setminus(X^1\cup X^2\cup W)$ with $(b-n-w)$ elements that satisfies $Z=\{z\in Z | L(\frac{b+c}{2})+1-z\}$.

For such $Z$, we have
\begin{equation}
\begin{aligned}
    \Delta(Z, X^2)&=\prod_{z\in Z, x_2\in X^2}|z-x_2|\\
    &=\prod_{z\in Z, x_1\in X^1}|(L(\frac{b+c}{2})+1-z)-(L(\frac{b+c}{2})+1-x_1)|\\
    &=\prod_{z\in Z, x_1\in X^1}|x_1-z|\\
    &=\Delta(Z, X^1)
\end{aligned}
\end{equation}
Similarly, we also have $\Delta(Z, Y^2)=\Delta(Z, Y^1)$.

Hence we have
\begin{equation}
\begin{aligned}
    \Delta(Z, X^1)=\sqrt{\Delta(Z, X^1)\Delta(Z, X^2)}&=\sqrt{\Delta(Z, X^1\cup X^2)}\\
    &=\sqrt{\Delta(Z, Y^1\cup Y^2)}\\
    &=\Delta(Z, Y^1)
\end{aligned}
\end{equation}

By same reasoning, we have $\Delta(X^1, W)=\Delta(Y^1, W)$.

Now, we observe a ratio $\frac{\Delta(Z\cup X^1\cup W \cup C)}{\Delta(Z\cup Y^1\cup W \cup C)}$ for any set $Z$:
\begin{equation}
\begin{aligned}
    &\frac{\Delta(Z\cup X^1\cup W \cup C)}{\Delta(Z\cup Y^1\cup W \cup C)}\\
    &=\frac{\Delta(Z)\Delta(X^1)\Delta(W)\Delta(C)}{\Delta(Z)\Delta(Y^1)\Delta(W)\Delta(C)}\\
    &\cdot\frac{\Delta(Z,X^1)\Delta(Z,W)\Delta(Z,C)\Delta(X^1,W)\Delta(X^1,C)\Delta(W,C)}{\Delta(Z,Y^1)\Delta(Z,W)\Delta(Z,C)\Delta(Y^1,W)\Delta(Y^1,C)\Delta(W,C)}\\
    &=\frac{s(X^1)}{s(Y^1)}\cdot\frac{\Delta(X^1,C)}{\Delta(Y^1,C)}
\end{aligned}
\end{equation}

Since this ratio does not depend on a choice of a set $Z$, by (3.35), (3.36) and (3.40), we have
\begin{equation}
\begin{aligned}
    &\frac{M_\odot(H_{a,b,c}^{0,\frac{b+c}{2}}(F(A_1:W_1),..., F(A_t:W_t) : m_1, m_2,..., m_t, m_{t+1}))}{M_\odot(H_{a,b,c}^{0,\frac{b+c}{2}}(B(A_1:W_1),..., B(A_t:W_t) : m_1, m_2,..., m_t, m_{t+1}))}\\
    &=\frac{s(X^1)}{s(Y^1)}\cdot\frac{\Delta(X^1,C)}{\Delta(Y^1,C)}
\end{aligned}
\end{equation}

However, as we have seen in the proof of the Theorem 2.3.,
\begin{equation}
\begin{aligned}
    &\frac{M(H_{a,b,c}^{0,\frac{b+c}{2}}(F(A_1:W_1),..., F(A_t:W_t) : m_1, m_2,..., m_t, m_{t+1}))}{M(H_{a,b,c}^{0,\frac{b+c}{2}}(B(A_1:W_1),..., B(A_t:W_t) : m_1, m_2,..., m_t, m_{t+1}))}\\
    &=\frac{s(X^1)s(X^2)}{s(Y^1)s(Y^2)}\cdot\frac{\Delta(X^1,C)}{\Delta(Y^1,C)}\cdot\frac{\Delta(B,X^2)}{\Delta(B,Y^2)}
\end{aligned}
\end{equation}

Since our region is centrally symmetric, we have
\begin{equation}
\begin{aligned}
    s(X^1)&=\frac{1}{H(p)}\prod_{x<y, x,y\in X^1}(y-x)\\
    &=\frac{1}{H(n)}\prod_{x<y, x,y\in X^1}((L(\frac{b+c}{2})+1-x)-(L(\frac{b+c}{2})+1-y))\\
    &=\frac{1}{H(n)}\prod_{y'<x', x',y'\in X^2}(x'-y')\\
    &=s(X^2)
\end{aligned}
\end{equation}
Similarly, $s(Y^1)=s(Y^2)$, $\Delta(B,X^2)=\Delta(X^1,C)$ and $\Delta(B,Y^2)=\Delta(Y^1,C)$

Hence we have
\begin{equation}
\begin{aligned}
    &\frac{M_\odot(H_{a,b,c}^{0,\frac{b+c}{2}}(F(A_1:W_1),..., F(A_t:W_t) : m_1, m_2,..., m_t, m_{t+1}))}{M_\odot(H_{a,b,c}^{0,\frac{b+c}{2}}(B(A_1:W_1),..., B(A_t:W_t) : m_1, m_2,..., m_t, m_{t+1}))}\\
    &=\frac{s(X^1)}{s(Y^1)}\cdot\frac{\Delta(X^1,C)}{\Delta(Y^1,C)}\\
    &=\sqrt{\frac{s(X^1)s(X^2)}{s(Y^1)s(Y^2)}\cdot\frac{\Delta(X^1,C)}{\Delta(Y^1,C)}\cdot\frac{\Delta(B,X^2)}{\Delta(B,Y^2)}}\\
    &=\sqrt{\frac{M(H_{a,b,c}^{0,\frac{b+c}{2}}(F(A_1:W_1),..., F(A_t:W_t) : m_1, m_2,..., m_t, m_{t+1}))}{M(H_{a,b,c}^{0,\frac{b+c}{2}}(B(A_1:W_1),..., B(A_t:W_t) : m_1, m_2,..., m_t, m_{t+1}))}}
\end{aligned}
\end{equation}

Also, by (3.21),
\begin{equation}
\begin{aligned}
    \frac{\Delta(X^1,C)}{\Delta(Y^1,C)}=&\prod_{i=1}^{t}\Bigg[\frac{H(d_{SE}(L^i_{u_i})H(d_{NE}(T^i)}{H(d_{SE}(T^i)H(d_{NE}(L^i_{u_i})}\\
    &\cdot\prod_{j<u_i, j\in J_i}\frac{H(d_{SE}(L^i_j)H(d_{NE}(R^i_j)}{H(d_{SE}(R^i_j)H(d_{NE}(L^i_j)}\prod_{j\geq u_i, j\in I_i}\frac{H(d_{SE}(R^i_j)H(d_{NE}(L^i_j)}{H(d_{SE}(L^i_j)H(d_{NE}(R^i_j)}\Bigg]
\end{aligned}
\end{equation}

Hence, by (3.43) and (3.44), we have
\begin{equation}
\begin{aligned}
    &\frac{M_\odot(H_{a,b,c}^{0,\frac{b+c}{2}}(F(A_1:W_1),..., F(A_t:W_t) : m_1, m_2,..., m_t, m_{t+1}))}{M_\odot(H_{a,b,c}^{0,\frac{b+c}{2}}(B(A_1:W_1),..., B(A_t:W_t) : m_1, m_2,..., m_t, m_{t+1}))}\\
    &=\sqrt{\frac{M(H_{a,b,c}^{0,\frac{b+c}{2}}(F(A_1:W_1),..., F(A_t:W_t) : m_1, m_2,..., m_t, m_{t+1}))}{M(H_{a,b,c}^{0,\frac{b+c}{2}}(B(A_1:W_1),..., B(A_t:W_t) : m_1, m_2,..., m_t, m_{t+1}))}}\\
    &=\frac{s(X^1)}{s(Y^1)}\cdot\prod_{i=1}^{t}\Bigg[\frac{H(d_{SE}(L^i_{u_i}))H(d_{NE}(T^i))}{H(d_{SE}(T^i))H(d_{NE}(L^i_{u_i}))}\\
    &\cdot\prod_{j<u_i, j\in J_i}\frac{H(d_{SE}(L^i_j))H(d_{NE}(R^i_j))}{H(d_{SE}(R^i_j))H(d_{NE}(L^i_j))}\prod_{j\geq u_i, j\in I_i}\frac{H(d_{SE}(R^i_j))H(d_{NE}(L^i_j))}{H(d_{SE}(L^i_j))H(d_{NE}(R^i_j))}\Bigg]
\end{aligned}
\end{equation}

\section{Acknowledgement}
The author would like to thank to his advisor, Professor Mihai Ciucu for his encouragement and useful discussions. The geometric interpretation of terms in formulas which unifies the results is due to him. Also, the author thanks Jeff Taylor for installing software and frequent helpful assistance.


\begin{thebibliography}{999}

\bibitem{C1}
  M. Ciucu,
  \emph{Enumeration of perfect matchings in graphs with reflective symmetry},
  J. Combin. Theory Ser. A \textbf{77} (1997),
  67-97

\bibitem{C2}
  M. Ciucu,
  \emph{Enumeration of lozenge tilings of punctured hexagons},
  J. Combin. Theory Ser. A \textbf{83} (1998),
  268-272

\bibitem{C3}
  M. Ciucu,
  \emph{The other dual of MacMahon's theorem on plane partitions},
  Advances in Mathematics \textbf{306} (2017),
  427-450

\bibitem{C4}
  M. Ciucu,
  \emph{Symmetries of shamrocks \RomanNumeralCaps{4}: The self-complementary case},
  preprint, arXiv:1906.02022

\bibitem{C5}
  M. Ciucu,
  \emph{Centrally symmetric tilings of fern-cored hexagons},
  preprint, arXiv:1906.02951

\bibitem{CEKZ}
  M. Ciucu, T. Eisenkolbl, C. Krattenthaler and D. Zare,
  \emph{Enumeration of lozenge tilings of hexagons with a central triangular hole},
  J. Combin. Theory A \textbf{95} (2001),
  251-334

\bibitem{CL}
  M. Ciucu, T. Lai,
  \emph{Lozenge Tilings of Doubly-intruded Hexagons},
  accepted to appear in J. Combin. Theory A

\bibitem{CJP}
  H. Cohn, M. Larsen and J. Propp,
  \emph{The shape of a typical boxed plane partition},
  New York J. of Math. \textbf{4} (1998),
  137-165

\bibitem{GT}
  I. M. Gelfand and M. L. Tsetlin,
  \emph{Finite-dimensional representation of the group of unimodular matrices (in Russian)},
  Doklady Akad. Nauk. SSSR (N. S.) \textbf{71} (1950),
  825-828

\bibitem{OK}
  S. Okada and C. Krattenthaler,
  \emph{The number of rhombus tilings of a "punctured" hexagon and the minor summations formula},
  Adv. Appl. Math. \textbf{21} (1998).
  381-404
  
\bibitem{R}
  H. Rosengren,
  \emph{Selberg integrals, Askey-Wilson polynomials and lozenge titlings of a hexagon with a triangular hole},
  J. Combin. Theory Ser. A \textbf{138} (2016),
  29-59
\end{thebibliography}
\end{document}